\documentclass[10pt]{amsart}
\usepackage{amssymb,amsbsy,amsmath,amsfonts,times}
\usepackage{latexsym,euscript,exscale}

\usepackage{helvet}

\title{On isometry groups and maximal symmetry}
\author{Valentin Ferenczi  and Christian Rosendal }
\address{Instituto de Matem\'atica e Estat\'istica \\
 Universidade de S\~ao Paulo \\
rua do Mat\~ao 1010 \\
Cidade Universit\'aria \\
05508-90 S\~ao Paulo, SP \\
Brazil.}
\email{ferenczi@ime.usp.br}
\address{Department of Mathematics, Statistics, and Computer Science (M/C 249)\\
University of Illinois at Chicago\\
851 S. Morgan St.\\
Chicago, IL 60607-7045\\
USA}
\email{rosendal.math@gmail.com}
\urladdr{http://www.math.uic.edu/$~$rosendal}

\date{}
\newcommand {\K}{\mathbb K}
\newcommand{\F}{\mathbb F}
\newcommand {\N}{\mathbb N}

\newcommand {\R}{\mathbb R}
\newcommand {\Z}{\mathbb Z}
\newcommand {\C}{\mathbb C}

\newcommand {\T}{\mathbb T}

\newcommand{\SSS}{\ku S}
\newcommand{\In}{{\ku I\!\!}n}
\newcommand{\AF}{\ku {A\!\! F}}

\newcommand{\norm}[1]{\lVert#1\rVert}
\newcommand{\Norm}[1]{\big\lVert#1\big\rVert}

\newcommand{\triple}[1]{|\!|\!|#1|\!|\!|}

\newcommand{\Id}{{\rm Id}}

\newcommand{\eps}{\epsilon}

\newcommand{\iso}{\cong}

\newcommand{\im}{{\rm im}}

\newcommand{\tom} {\emptyset}

\newcommand{\saa}{\Rightarrow}

\newcommand{\til}{\rightarrow}

\newcommand{\Lim}[1]{\mathop{\longrightarrow}\limits_{#1}}

\newcommand {\Del}{ \; \Big| \;}
\newcommand {\del}{ \; \big| \;}

\newcommand {\go} {\mathfrak}
\newcommand {\ku} {\mathcal}

\newcommand{\inv}{^{-1}}

\newcommand {\e} {\exists}
\renewcommand {\a} {\forall}

\newtheorem{thm}{Theorem}[section]
\newtheorem{cor}[thm]{Corollary}
\newtheorem{lemme}[thm]{Lemma}
\newtheorem{prop} [thm] {Proposition}
\newtheorem{defi} [thm] {Definition}

\newtheorem{prob}[thm]{Problem}

\newtheorem{claim}[thm] {Claim}

\newtheorem{conj}[thm]{Conjecture}

\usepackage{fouriernc}

\thanks{V. Ferenczi acknowledges the support of CNPq, project 452068/2010-0, and FAPESP, projects 2010/05182-1 and 2010/17493-1. C. Rosendal was supported by NSF grant DMS 0901405}

\begin{document}

\subjclass[2000]{Primary: 22F50, 46B03, 46B04. Secondary: 03E15}

\keywords{Isometries of Banach spaces, maximal norms, almost transitive norms}

\maketitle

\begin{abstract}
We study problems of maximal symmetry in Banach spaces. This is done by providing an analysis of the structure of small subgroups of the general linear group $GL(X)$, where $X$ is a separable reflexive Banach space. In particular, we provide the first known example of a Banach space $X$ without any equivalent maximal norm, or equivalently such that $GL(X)$ contains no maximal bounded subgroup. Moreover, this space $X$ may be chosen to be super-reflexive.
\end{abstract}

\

\begin{center}
In memory of 
Greg Hjorth (1963--2011)\\
and Nigel Kalton (1946--2010).
\end{center}

\

\tableofcontents


\section{Introduction}

\subsection{Maximal norms and the problems of Mazur and Dixmier}
Two outstanding problems of functional analysis are S. Mazur's {\em  rotation problem}, asking whether any separable Banach space whose isometry group acts transitively on the sphere must be Hilbert space, and J. Dixmier's {\em  unitarisability problem}, asking whether any countable group, all of whose bounded representations on Hilbert space are unitarisable, must be amenable. 
Though these problems do not on the surface seem to be related, they both point to common geometric aspects of Hilbert space that are far from being well-understood.

Mazur's problem, which can be found in S. Banach's classical work \cite{Banach} as well as in the Scottish Book \cite{mauldin}, is perhaps best understood as two separate problems, both of which remain open to this day.

\begin{prob}[Mazur's rotation problem, first part]
Suppose $X$ is a separable Banach space whose isometry group acts transitively on the sphere $S_X$. Is $X$ Hilbertian, i.e., isomorphic to the separable Hilbert space $\ku H$?
\end{prob}

Remark that, in order for the problem to be non-trivial, the separability condition in the hypothesis is necessary. For the isometry group of $L_p$ induces a dense orbit on the sphere  (see \cite{R,GJK}) and thus the isometry group will act transitively on any ultrapower of $L_p$, which itself is an $L_p$-space. 
\begin{prob}[Mazur's rotation problem, second part]
Suppose $\triple\cdot$ is an equivalent norm on $\ku H$ such that ${\rm Isom}(\ku H,\triple\cdot)$ acts transitively on the unit sphere 
$S_{\ku H}^{\triple\cdot}$. Is $\triple\cdot$ necessarily euclidean?
\end{prob}
Establishing a vocabulary to study these problems, A. Pe\l czy\'nski and S. Rolewicz \cite{PR} (see also Rolewicz' book \cite{R}) defined the norm $\norm\cdot$ on a Banach space $X$ to be {\em maximal} if whenever $\triple\cdot$ is an equivalent norm on $X$ with 
$$
{\rm Isom}(X,\norm\cdot)\leqslant {\rm Isom}(X,\triple\cdot),
$$ 
then
$$
{\rm Isom}(X,\norm\cdot)={\rm Isom}(X,\triple\cdot).
$$ 
Also, the norm is {\em transitive} if the isometry group acts transitively on the unit sphere.

Thus, a norm is maximal if one cannot replace it by another equivalent norm that has strictly more isometries, or, more suggestively, if the unit ball $B_X^{\norm\cdot}$ is a maximally symmetric body in $X$. Note also that if $\norm \cdot$ is transitive, then $S_X^{\norm\cdot}$ is the orbit of a single point $x\in S_X^{\norm\cdot}$ under the action of ${\rm Isom}(X,\norm\cdot)$, and so any proper supergroup $G$ of ${\rm Isom}(X,\norm\cdot)$ in $GL(X)$ must send $x$ to some point $\lambda x$ for $|\lambda|\neq 1$, from which it follows that $G$ cannot be a group of isometries for any norm. So transitivity implies maximality and thus the standard euclidean norm $\norm\cdot_2$ is maximal on $\ku H$. Also, Mazur \cite{mazur} showed himself that any transitive norm on a finite-dimensional space is, in effect, euclidean (see the survey papers  by F. Cabello-S\'anchez \cite{CS} and J. Becerra Guerrero and A. Rodr\'iguez-Palacios \cite{becerra} for more information on the rotation problem and maximal norms).

Another way of understanding these concepts, pointing towards the unitarisability problem of Dixmier, is by considering the $G$-invariant norms corresponding to a bounded subgroup $G\leqslant GL(X)$. Here $G$ is {\em bounded} if $\norm G=\sup_{T\in G}\norm T<\infty$. Note first that if $G$ is bounded, then
$$
\triple x=\sup_{T\in G}\norm{Tx}
$$
defines an equivalent $G$-invariant norm on $X$, i.e., $G\leqslant {\rm Isom}(X, \triple\cdot)$. Moreover, if $\norm\cdot$ is uniformly convex, then so is $\triple\cdot$ (see, e.g., Proposition 2.3 in \cite{furman}). However, if $X=\ku H$ and $\norm\cdot$ is euclidean, i.e., induced by an inner product, then $\triple\cdot$ will not, in general, be euclidean. The question of which bounded subgroups of $GL(\ku H)$ admit invariant euclidean norms has a long history. It is a classical result of representation theory dating back to the beginning of the 20th century that if $G\leqslant GL(\C^n)$ is a bounded subgroup, then there is a $G$-invariant inner product, thus inducing a $G$-invariant euclidean norm. Also, in the 1930s, B. Sz.-Nagy showed that any bounded representation $\pi\colon \Z\til GL(\ku H)$ is {\em unitarisable}, i.e., $\ku H$ admits an equivalent $\pi(\Z)$-invariant inner product and, with the advent of amenability in the 1940s,  this was extended by M. Day \cite{day} and J. Dixmier \cite{dixmier} to any bounded representation of an amenable topological group via averaging over an invariant mean. In the opposite direction, L. Ehrenpreis and F. I. Mautner \cite{mautner} constructed a non-unitarisable bounded representation of $SL_2(\R)$ on $\ku H$ and the group $SL_2(\R)$ was later replaced by any countable group containing the free group $\F_2$. However, since $\F_2$ does not embed into all non-amenable countable groups, the question of whether the  result of Sz.-Nagy, Day and Dixmier reverses still remains pertinent.

\begin{prob}[Dixmier's unitarisability problem] 
Suppose $\Gamma$ is a countable group all of whose bounded representations on $\ku H$ are unitarisable. Is $\Gamma$ amenable?
\end{prob}

Nevertheless, though not every bounded representation of $\F_2$ in $GL(\ku H)$ is unitarisable, it still seems to be unknown whether all of its bounded representations admit equivalent invariant maximal or even transitive norms. And similarly, while the second part of the rotation problem asks whether any equivalent transitive norm on $\ku H$ is euclidean, the stronger question of whether any equivalent maximal norm is euclidean also remains open. Of course, the counter-example of Ehrenpreis and Mautner limits how much of these two questions can hold simultaneously (see the recent paper  \cite{pisier} by G. Pisier for material on the current status of Dixmier's problem).

In a more general direction, the work of Pe\l czy\'n\-ski and Rolewicz led people to investigate which spaces have maximal norms. Since any bounded subgroup $G\leqslant GL(X)$ is a group of isometries for an equivalent maximal norm, one observes that a norm $\norm\cdot$ is maximal if and only if the corresponding isometry group is a maximal bounded subgroup of $GL(X)$. Thus, in analogy to the existence of maximal compact subgroups of semi-simple Lie groups, it is natural to suspect that a judicious choice of smoothing procedures on a space $X$ could eventually lead to a most symmetric norm, which then would be maximal on $X$. 
But, even so, fundamental questions on maximal norms have remained open, including notably the longstanding problem, formulated by G. Wood in \cite{W}, whether any Banach space $X$ admits an equivalent maximal norm. In fact, even the question of whether any bounded $G\leqslant GL(X)$ is contained in a maximal bounded subgroup was hitherto left unresolved. 
Our main result answers these, as well as another problem of R. Deville, G. Godefroy and V. Zizler \cite{DGZ}, in the negative.
\begin{thm}\label{main1}
There is a separable super-reflexive Banach space $X$ such that $GL(X)$ contains no maximal bounded subgroups, i.e., $X$ has no equivalent maximal norm.
\end{thm}


\subsection{Non-trivial isometries of Banach spaces}
A second motivation and a source of tools for our work comes from the seminal construction of W. T. Gowers and B. Maurey \cite{GM} of a space $GM$ with a small algebra of operators, namely, such that any operator on $GM$ is a strictly singular perturbation of a scalar multiple of the identity map. The currently strongest result in this direction, due to S. A. Argyros and R. G. Haydon \cite{AH}, is the construction of a Banach space $AH$ on which every operator is a compact perturbation of a scalar multiple of the identity. Furthermore, since AH has a Schauder basis, every compact operator is a limit of operators of finite rank.

These results largely answer the question of whether every Banach space admits non-trivial operators, but one can ask the same question for isometries, i.e., does every Banach space admit a non-trivial surjective isometry? After partial answers by P. Semenev and A. Skorik \cite{SS} and S. Bellenot \cite{Bel}, one version of this question was answered in the negative by K. Jarosz \cite{J}. Jarosz proved that any real or complex Banach space admits an equivalent  norm with only trivial isometries, namely, such that any surjective isometry is a scalar multiple of the identity, $\lambda\Id$, for $|\lambda|=1$. Thus, no isomorphic property of a space can force the existence of a non-trivial surjective linear isometry.

Of course, this does not prevent the group of isometries to be extremely non-trivial in some other equivalent norm. So one would like results relating the size of the isometry group ${\rm Isom}(X,\norm\cdot)$ with the isomorphic structure of $X$.
Let us first remark that any infinite-dimensional Banach space $X$ can always be equivalently renormed such that $X=F \oplus_1 H$, where $F$ is a finite-dimensional euclidean space. So, in this case, ${\rm Isom}(X)$ will at least contain a subgroup isomorphic to ${\rm Isom}(F)$. Actually, if $X$ is a separable, infinite-dimensional, real space and $G$ is a finite group, then it is possible to find an equivalent norm for which $\{-1,1\} \times G$ is isomorphic to the group of isometries on $X$  \cite{FG}.

Thus, allowing for renormings, we need a less restrictive concept of when an isometry is trivial.
\begin{defi} 
A bounded  subgroup $G\leqslant GL(X)$ acts {\em nearly trivially} on $X$ if there is a $G$-invariant decomposition  $X=F \oplus H$, where $F$ is finite-dimensional and $G$ acts by trivial isometries on $H$.
\end{defi}

As an initial step towards Theorem \ref{main1}, we show that in a certain class of spaces, each individual isometry acts nearly trivially. For that, we shall need to improve on some earlier work of to F. R\"abiger and W. J. Ricker \cite{RR,RR2}. By their results, any isometry of a so called hereditarily indecomposable complex Banach space is a compact perturbation of a scalar multiple of the identity, but this can be improved as follows.

\begin{thm}\label{nearly trivial-intro}
Let $X$ be a Banach space containing no unconditional basic sequence and on which every operator is of the form $\lambda\Id+S$, for $S$ strictly singular. Then each individual isometry acts nearly trivially on $X$.
\end{thm}

The main problem is then to investigate when we can proceed from single isometries acting nearly trivially to an understanding of the global structure of the isometry group ${\rm Isom}(X)$. Disregarding for the moment the scalar multiples of the identity on $X$, we consider the automorphisms of $X$ that individually acts nearly trivially on $X$. For this, we let $GL_f(X)$ denote the subgroup of $GL(X)$ consisting of all automorphisms of the form
$$
\Id+A,
$$
where $A$ is a finite-rank operator on $X$. We then establish, in the case of separable reflexive $X$, the structure of  bounded subgroups of $GL_f(X)$ that are strongly closed in $GL(X)$. The following statement refers to the strong operator topology on $G$.

\begin{thm}\label{decomp-intro}
Suppose $X$ is a separable, reflexive Banach space and $G\leqslant GL_f(X)$ is bounded and strongly closed in $GL(X)$. 
Let also $G_0\leqslant G$ denote the connected component of the identity in $G$.

Then $G_0$ acts nearly trivially on $X$ and therefore is a compact Lie group. Moreover, $G_0$ is open in $G$, while $G/G_0$ is a countable, locally finite group. It follows that $G$ is an amenable Lie group.

Furthermore, $X$ admits a $G$-invariant decomposition $X=X_1\oplus X_2\oplus X_3\oplus X_4$, where
\begin{enumerate}
\item no non-zero point of $X_1$ has a relatively compact $G$-orbit,
\item every $G$-orbit on $X_2\oplus X_3$ is relatively compact,
\item $X_4$ is the subspace of points which are fixed by $G$,
\item $X_2$ is finite-dimensional and $X_1\oplus X_3\oplus X_4$ is the subspace of points which are fixed by $G_0$,
\item if $X_1\neq\{0\}$, then $X_1$ has a complemented subspace with a Schauder basis, while if $X_1=\{0\}$, then  $G$ acts nearly trivially on $X$.
\end{enumerate}
\end{thm}

Combining Theorem \ref{nearly trivial-intro} with a somewhat simpler version of Theorem  \ref{decomp-intro} along with an earlier construction of a super-reflexive hereditarily indecomposable space due to the first named author \cite{F}, we are able to conclude the following result from which Theorem \ref{main1} is easily obtained.

\begin{thm}\label{main2} 
Let $X$ be a separable, reflexive, hereditarily indecomposable, complex Banach space without a Schauder basis. Then for any equivalent norm on $X$, the group of isometries acts nearly trivially on $X$.

Moreover, there are super-reflexive spaces satisfying these hypotheses.
\end{thm}


\section{Notation and complexifications}

\subsection{Some notation and terminology}

For a  Banach space $X$, we denote by ${\mathcal L}(X)$ the algebra of continuous linear operators on $X$,  by $GL(X)$ the {\em general linear group} of $X$, i.e., the group of all continuous linear automorphisms of $X$, and by ${\rm Isom}(X)$ the group of surjective linear isometries of $X$. In the remainder of the paper, unless explicitly stated otherwise, an isometry of $X$ is always assumed to be surjective, so we shall not state this hypothesis explicitly. We also denote the unit sphere of $X$ by $S_X$ and the closed unit ball by $B_X$.

If $X$ is a complex space and $T$ is an operator on $X$, $\sigma(T)$ denotes the spectrum of $T$ and one  observes that $\sigma(T)$ is a subset of the unit circle $\T$ whenever $T$ is an isometry.
On the other hand, if $T$ is compact, then $T$ is a Riesz operator, which means that $\sigma(T)$ is either a finite sequence of eigenvalues with finite multiplicity together with $0$, or an infinite converging sequence of such eigenvalues together with the limit point $0$.
 
As our proofs use methods from representation theory, spectral theory, renorming theory, as well as general Banach space theory, we have tried to give self-contained and detailed proofs of our results in order for the paper to remain readable for a larger audience. Of course some of the background material is true in a broader setting, which may be easily found by consulting the literature. Our main references will be the book of N. Dunford and J. Schwarz \cite{DS} for spectral theory, the book of R. Deville, G. Godefroy and V. Zizler \cite{DGZ} for renorming theory, and the books of J. Lindenstrauss and L. Tzafriri \cite{LT} and of Y. Benyamini and J. Lindenstrauss \cite{BL} for general Banach space and some operator theory. We also recommand the book of R. Fleming and J. Jamison \cite{FJ}, especially Chapter 12, for more information on isometries of Banach spaces.

\subsection{Complex spaces versus real spaces}\label{complexification}

Our main results will be valid both in the real and the complex settings, though different techniques will sometimes be needed to cover each separate case.

For the part of our demonstrations using spectral theory, as is classical, we shall first prove our results in the complex case and thereafter use complexification to extend them to the real case. We recall briefly how this is done and the links that exist between a real space and its complexification.

If $X$ is a real Banach space with norm $\norm\cdot$, the complexification $\hat{X}$ of $X$ is defined as the Cartesian square $X \times X$, whose elements are written $x+iy$ rather than $(x,y)$ for $x,y \in X$, equipped with the complex scalar multiplication given by
$$
(a+ib)\cdot(x+iy)=ax-by+i(bx+ay),
$$ 
for $a,b\in \R$ and $x,y\in X$, and with the equivalent norm
$$
\triple{x+iy}=\sup_{\theta \in [0,2\pi]} \|e^{i\theta}(x+iy)\|_2,
$$
where
$$\|x+iy\|_2=\sqrt{\norm{x}^2+\norm{y}^2}.$$
Any operator in ${\mathcal L}(\hat{X})$ may be written of the form $T+iU$, where $T, U$ belong
to ${\mathcal L}(X)$, that is, for $x,y \in X$,
$$(T+iU)(x+iy)=Tx-Uy+i(Ux+Ty).$$
We denote by $c$ the natural isometric homomorphism from ${\mathcal L}(X)$ into ${\mathcal L}(\hat{X})$ associating to $T$ the operator $\hat{T}$ defined by
$$\hat{T}=T+i0.$$

It is then straightforward to check that the image by $c$ of an automorphism (respectively isometry, finite-rank perturbation of the identity, compact perturbation of the identity) of $X$ is an automorphism (respectively isometry, finite-rank perturbation of the identity, compact perturbation of the identity) of $\hat{X}$. 
In other words,  the map $c$ provides an embedding of natural subgroups of $GL(X)$ into their counterparts in $GL(\hat{X})$.

Conversely, in renorming theory it is usually assumed that the spaces are real. One obtains renorming results on a complex space $X$ simply by considering only its $\R$-linear structure which defines a real space $X_{\R}$. It is well-known that many isomorphic properties of a space do not depend on it being seen as real or complex, and we recall here briefly the facts that we shall rely on later.

First, by \cite{BL}, a space $X$ has the Radon--Nikodym Property if and only if every Lipschitz function from $\R$ into $X$ is differentiable almost everywhere, which only depends on the $\R$-linear structure of $X$. A space $X$ is reflexive if and only if the closed unit ball of $X$ is weakly sequentially compact, and hence $X$ is reflexive if and only if $X_\R$ is reflexive.

Moreover, the map $\phi \mapsto {\rm Re}(\phi)$ is  an $\R$-linear isometry from
$X^*$ onto $(X_\R)^*$ with inverse $\psi\mapsto \phi$ defined by
$$
\phi(x)=\psi(x)-i\psi(ix).
$$
So the dual norms on $X^*$ and $(X_\R)^*$ coincide up to this identification and  $X$ has separable dual if and only if $X_\R$ has. Likewise, when  $\phi$ is a $\C$-support functional for $x_0 \in   S_X$,  ${\rm Re}(\phi)$ is an $\R$-support functional for $x_0$. Thus, the above identification shows that if in a complex space a point $x$ in $S_X$ has a unique $\R$-support functional, then it has a unique $\C$-support functional, which, in particular,  will happen when the norm on $X$ is G\^ateaux-differentiable \cite{DGZ}.


\section{Bounded subgroups of $GL(X)$}

In the present section, we shall review some general facts about bounded subgroups of $GL(X)$, so, apart from Theorems \ref{in then af} and \ref{christian},  all of the material here is well-known, but maybe hard to find in any single source.

\subsection{Topologies on $GL(X)$}

Suppose $X$ is a real or complex Banach space and $G\leqslant GL(X)$ is a weakly bounded subgroup, i.e., such that
for any $x\in X$ and $\phi\in X^*$,
$$
\sup_{T\in G}|\phi(Tx)|<\infty.
$$
Then, by the uniform boundedness principle, $G$ is actually norm bounded, that is,
$$
\norm{G}=\sup_{T\in G}\|T\|<\infty.
$$
So without ambiguity we can simply refer to $G$ as a {\em bounded} subgroup of $GL(X)$.

Note that if $G$ is bounded, then
$$
\triple{x}=\sup_{T\in G}\|Tx\|
$$
is an equivalent norm on $X$ such that $G$ acts by {\em isometries} on $(X,\triple{\cdot})$.
Therefore, bounded subgroups of $GL(X)$ are simply groups of isometries for equivalent norms on $X$.

Let us also stress the fact that, although the operator norm changes when $X$ is given an equivalent norm, the norm,  weak and strong operator topologies on $GL(X)$ remain unaltered. 

Recall that if $X$ is separable, the isometry group, ${\rm Isom}(X)$, is a {\em Polish} group in the strong operator topology, i.e., a separable topological group whose topology can be induced by a complete metric. Since any strongly closed bounded subgroup $G\leqslant GL(X)$ can be seen as a strongly closed subgroup of ${\rm Isom}(X)$ for an equivalent norm on $X$, provided $X$ is separable, we find that $G$ is a closed subgroup of a Polish group and hence is Polish itself.

Note also that the norm  induces an invariant, complete metric on ${\rm Isom}(X)$, that is, $\norm {TSU-TRU}=\norm{S-R}$ for all $T,S,R,U\in {\rm Isom}(X)$, and so ${\rm Isom}(X)$ and, similarly, any bounded subgroup $G\leqslant GL(X)$, is a {\em SIN} group in the norm topology, i.e., admits a neighbourhood basis at the identity consisting of conjugacy invariant sets. 

Of course, even when $X$ is separable, the norm topology can be non-separable on ${\rm Isom}(X)$, but, as we shall see, for certain small subgroups of $GL(X)$ it coincides with the strong operator topology, which allows for an interesting combination of different techniques.

Note that if $X$ is separable, we can choose a dense subset of the unit sphere $\{x_n\}\subseteq   S_X$ and corresponding norming functionals $\{\phi_n\}\subseteq   S_{X^*}$, $\phi_n(x_n)=1$. Using these, we can write the closed unit ball $B_X$ as 
$$
B_X=\bigcap_{n,m}\{x\in X\del \phi_n(x)<1+ 1/m\},
$$
which shows that $B_X$ is a countable intersection of open half-spaces in $X$ and similarly for any other closed ball in $X$. Thus, if  $\pi\colon G\til GL(X)$ is a bounded, weakly continuous representation of a Polish group $G$, then for any $x\in X$ and $\eps>0$, the set
$$
\{g\in G\del \pi(g)x\in  \overline{B(x,\eps)}\}
$$
is a countable intersection of open sets in $G$ and hence is Borel. It follows that $\pi$ is a Borel homomorphism from a Polish group into the separable group $(\pi(G), {\rm SOT})$ and hence, by Pettis's Theorem (see \cite{kechris}, (9.10)), $\pi$ is strongly continuous.

\begin{prop}\label{weakstrong}
Let $X$ be a separable Banach space and let $\pi\colon G\til GL(X)$ be a bounded, weakly continuous representation of a Polish group $G$. Then $\pi$ is strongly continuous.
\end{prop}

In fact, if $X^*$ is separable, the weak and strong operator topologies coincide on any bounded subgroup $G\leqslant GL(X)$ (see, e.g., \cite{megrelishvili} for more information on this).

What is more important is that if $\pi\colon G\til GL(X)$ is a strongly continuous bounded representation of a Polish group, the induced dual representation $\pi^*\colon G\til GL(X^*)$ is in general only {\em ultraweakly} continuous, i.e., $\pi^*(g_i)(\phi)\Lim{{\rm weak}^*}\pi^*(g)(\phi)$ for $\phi\in X^*$ and $g_i\til g$. Of course, if $X$ is separable, reflexive, this means that $\pi^*$ is weakly continuous and thus also strongly continuous.

In the following, the default topology on $GL(X)$ and its subgroups is the strong operator topology. So, unless otherwise stated, all statements refer to this topology. Moreover, if $G\leqslant GL(X)$, then $\overline{G}^{\rm SOT}$ refers to the strong closure in $GL(X)$ and {\em not} in $\ku L(X)$. This is important, since even a bounded subgroup that is strongly closed in $GL(X)$ may not be strongly closed in $\ku L(X)$, e.g., the unitary group of infinite-dimensional Hilbert space, $U(\ell_2)$, is not strongly closed in $\ku L(\ell_2)$. 
This is in opposition to the well-known fact that any {\em bounded} subgroup that is norm closed in $GL(X)$ is also norm closed in $\ku L(X)$. However, for potentially unbounded $G\leqslant GL(X)$, we let $\overline{G}^{\|\cdot\|}$ denote the norm closure of $G$ in $GL(X)$.

A topological group $G$ is said to be {\em precompact} if any non-empty open set $\ku U\subseteq G$ covers $G$ by finitely many left translates, i.e., $G=A\ku U$ for some finite set $A\subseteq G$. For Polish groups, this is equivalent to being compact, but, e.g., any non-closed subgroup of a compact Polish group is only precompact and not compact.

For the next proposition, we recall that a vector $x\in X$ is said to be {\em almost periodic}, with respect to some $G\leqslant GL(X)$, if the $G$-orbit of $x$ is relatively compact or, equivalently, totally bounded in $X$.  In analogy with this, $G$ is said to be {\em almost periodic} if every $x\in X$ is almost periodic.

\begin{prop}\label{Gcompact}
Suppose $X$ is a  Banach space and $G\leqslant GL(X)$ is a bounded subgroup. Then the following are equivalent.
\begin{enumerate}
\item $X$ is the closed linear span of its finite-dimensional irreducible  subspaces,
\item $G$ is almost periodic,
\item $G$ is precompact,
\item $\overline{G}^{\rm SOT}$ is compact.
\end{enumerate}
\end{prop}

\begin{proof}(1)$\saa$(2): It is an easy exercise to see that the set of almost periodic points form a closed linear subspace and, moreover, since $G$ is bounded, any finite-dimensional $G$-invariant subspace is contained in the set of almost periodic points. So (2) follows from (1).

(2)$\saa$(3): Note that if $G$ is not precompact in the strong operator topology, then there is an open neighbourhood $\ku U$ of $\Id$ that does not cover $G$ by a finite number of left translates. It follows that we can find a finite sequence of normalised vectors $x_1,\ldots,x_n\in X$, $\eps>0$ and an infinite set $\ku A\subseteq G$ such that for distinct $T,U\in \ku A$ there is $i$ with $\|Tx_i-Ux_i\|>\eps$.
But then, by the infinite version of Ramsey's theorem \cite{ramsey}, we may find some $i$ and some infinite subset $\ku B$ of $\ku A$ such that $\|Tx_i-Ux_i\| > \eps$ whenever $T,U\in \ku B$ are distinct, which shows that the $G$-orbit of $x_i$ is not relatively compact.

(3)$\saa$(4): If $G$ is precompact, $\overline{G}^{\rm SOT}$ is easily seen to be precompact. It  follows that every $\overline{G}^{\rm SOT}$-orbit is totally bounded, i.e., relatively compact. So, to see that $\overline{G}^{\rm SOT}$ is compact, let $(g_i)$ be a net in $\overline{G}^{\rm SOT}$ and pick a subnet $(h_j)$ such that, for every $x\in X$, $(h_jx)$ and $(h_j\inv x)$ converge to some $Tx$ and $Sx$ respectively. It follows that $S=T\inv\in\overline{G}^{\rm SOT}$ and so $(h_j)$ converges in the strong operator topology to $T\in \overline{G}^{\rm SOT}$. Since every net has a convergent subnet, $\overline{G}^{\rm SOT}$ is compact.

(4)$\saa$(1): Suppose $\overline G^{\rm SOT}\leqslant GL(X)$ is compact and consider the tautological strongly continuous representation $\pi\colon \overline G^{\rm SOT}\til GL(X)$. By a result going back to at least K. Shiga \cite{shiga}, since $\overline G^{\rm SOT}$ is compact, $X$ is the closed linear span of its finite-dimensional irreducible subspaces, i.e., minimal non-trivial $G$-invariant subspaces and so (1) follows. (Note that the result of Shiga is stated only for the complex case in \cite{shiga}, but the real case follows from considering the complexification).
\end{proof}

At several occasions we shall be using the following theorem due to I. Gelfand (see \cite{HP}): If $T$ is an element of a  complex unital   Banach algebra $\go A$, e.g., $\go A=\ku L(X)$, with $\sigma(T)=\{1\}$ and $\sup_{n\in \Z}\|T^n\|<\infty$, then $T=1$.

Since $GL(X)$ is a norm open subset of the Banach space $\ku L(X)$, it is a {\em Banach--Lie} group, but is of course far from being a (finite-dimensional) Lie group. 
\begin{thm}\label{NSS}
Let $X$ be a Banach space. Then, in the norm topology, $GL(X)$ has {\em no small subgroups}, that is, there is $\eps>0$ (in fact $\eps=\sqrt2$) such that
$$
\{T\in GL(X)\del \|T-\Id\|<\eps\}
$$
contains no non-trivial subgroup.

If follows that if $G\leqslant GL(X)$ is locally compact, second countable in the norm-topology, then $G$ is a Lie group.
\end{thm}

\begin{proof}
Assume first that $X$ is a complex space.
We claim that for any $T \in GL(X)$ such that $(T^n)_{n \in \Z}$ is bounded, any $\lambda\in \sigma(T)$ is an approximate eigenvalue, that is, $Tx_n-\lambda x_n\til 0$ for some $x_n\in   S_X$. For otherwise, $T-\lambda\Id$ is bounded away from $0$ and hence will be an embedding of $X$ into $X$ whose range is a closed proper subspace of $X$.  Since $X$ can be renormed so that $T$ is an isometry, we have that $\sigma(T)\subseteq \T$, and so we may find find $\lambda_n\notin \sigma(T)$ such that  $\lambda_n\til \lambda$.  Therefore, if we choose $y\notin \im(T-\lambda\Id)$, there are $x_n\in X$ such that $Tx_n-\lambda_nx_n=y$ and so either $\|x_n\|$ is bounded or can be assumed to tend to infinity. In the second case, we see that for $z_n=\frac {x_n}{\|x_n\|}$ one has $Tz_n-\lambda_nz_n=\frac{y}{\|x_n\|}\til 0$ and so also $Tz_n-\lambda z_n\til 0$, contradicting that $T-\lambda\Id$ is bounded away from $0$. And, in the first case, 
$$
\|(Tx_n-\lambda x_n)-y\|=\|(Tx_n-\lambda x_n)-(Tx_n-\lambda_n x_n)\|=|\lambda_n-\lambda|\cdot\|x_n\|\til 0,
$$
contradicting that $y$ is not in the closed subspace $\im (T-\lambda\Id)$.

 Now let $T\in GL(X)$ satisfy $\|T^n-\Id\|<\sqrt2$ for all $n\in \Z$.
If $\lambda\in \sigma(T)$, then $\lambda^n$ is an approximate eigenvalue of $T^n$ for any $n \in \N$, and so it follows that $|\lambda^n-1|<\sqrt2$ for all $n\in \N$ and therefore that $\lambda=1$. So  $\sigma(T)=\{1\}$. It suffices now to apply Gelfand's Theorem to conclude that $T=\Id$. We have thus shown that $\big\{T\in GL(X)\del \|T-\Id\|<\sqrt2\big\}$ contains no non-trivial subgroup.

If $X$ is a real Banach space, it suffices again to consider the complexification of $X$.

For the second part of the theorem, we note that by the Gleason - Montgomery - Yamabe - Zippin solution to Hilbert's 5th problem (see, e.g.,  \cite{kaplansky} for an exposition), any locally compact, second countable group with no small subgroups is a Lie group.
\end{proof}

\subsection{Ideals and subgroups}

Note that when $\ku I\subseteq \ku L(X)$ is a two-sided operator ideal, the subgroup $GL_{\ku I}(X)\leqslant GL(X)$ consisting of all  $\ku I$-perturbations of the identity, that is, invertible operators of the form 
$$
T=\Id+A,
$$
where $A\in \ku I$, is normal in $GL(X)$. Moreover, if $\ku I$ is norm closed in $\ku L(X)$, then $GL_\ku I(X)$ is a norm closed subgroup of $GL(X)$.

Of particular importance for our investigation are the ideals
of respectively finite-rank, almost finite-rank, compact, strictly singular and inessential operators. Namely,
\begin{itemize}
\item $\ku F(X)=\{T\in \ku L(X)\del T \text{ has finite-rank}\}$,
\item ${\AF}(X)=\overline{\ku F(X)}^{\norm\cdot}$,
\item $\ku K(X)=\{T\in \ku L(X)\del T\text{ is compact}\}$,
\item ${\SSS}(X)=\{T\in \ku L(X)\del T \text{ is strictly singular}\}$,
\item ${\In}(X)=\{T\in\ku L(X)\del T\text{ is inessential}\}$.
\end{itemize}
Here an operator $T\in \ku L(X)$ is said to be {\em strictly singular} if there is no infinite-dimensional subspace $Y\subseteq X$ such that $T\colon Y\til X$ is an isomorphic embedding. Also, $T\in \ku L(X)$ is {\em inessential} if for any $S\in \ku L(X)$, the operator $\Id+ST$ is Fredholm, i.e., has closed image, finite-dimensional kernel and finite co-rank. In particular, for any $T\in {\In}(X)$ and $t\in [0,1]$, $\Id+tT$ is Fredholm and $\Id+T$ must have Fredholm index $0$, since the index is norm continuous and $(\Id+tT)_{t\in[0,1]}$ is a continuous path from $\Id$ to $\Id+T$.  More information about the ideal of inessential operators may be found in \cite{Gon}. We then have the following inclusions
$$
\ku F(X)\subseteq \AF(X)\subseteq \ku K(X)\subseteq \SSS(X)\subseteq \In(X),
$$
which gives us similar inclusions between the corresponding subgroups of $GL(X)$, that, for simplicity, we shall denote respectively by
$$
GL_f(X)\;\subseteq\; GL_{af}(X)\;\subseteq \;GL_c(X)\;\subseteq \;GL_{s}(X)\;\subseteq\; GL_{in}(X).
$$ 
We also note that the ideals $\AF(X)$, $ \ku K(X)$ and  $\SSS(X)$ are norm closed in $\ku L(X)$.  $GL_c(X)$ is usually called the {\em Fredholm group}, though sometimes this refers more specifically to $GL_c(\ell_2)$.

Now, since the compact operators form the only non-trivial norm closed ideal of ${\mathcal L}(\ell_2)$, we have that  
$$
GL_{af}(\ell_2)=GL_{c}(\ell_2)=GL_{s}(\ell_2).
$$
Similarly, if $X$ has the approximation property, then $\ku{K}(X)=\AF(X)$ and hence $GL_{af}(X)=GL_c(X)$. Though these equalities do not hold for general Banach spaces, as we shall see now, any {\em bounded} subgroup of $GL_{in}(X)$ is contained in $GL_{af}(X)$, so from our perspective, there is no loss of generality in only considering $GL_{af}(X)$. It should be noted that Theorem \ref{in then af} generalises and simplifies results of R\"abiger--Ricker  \cite{RR,RR2} and Ferenczi--Galego \cite{FG}.

\begin{thm}\label{in then af}
Let $X$ be a Banach space and $G\leqslant GL_{in}(X)$ be a bounded subgroup. Then $G$ is contained in $GL_{af}(X)$.
\end{thm}

\begin{proof} 
It suffices to show that if $T \in GL_{in}(X)$ and  $\{T^n \del n \in \Z\}$
is bounded, then $T \in {GL_{af}(X)}$. For this, we may assume that $X$ is infinite-dimensional.

Suppose first that $X$ is complex and work in the norm topology on $\ku L(X)$. Consider the quotient algebra
$$
{\mathcal B}=\ku L(X)/\AF(X),
$$
and let 
$$
\alpha\colon\ku L(X) \rightarrow {\mathcal B}
$$ 
be the corresponding quotient map.

Fix $T=\Id+U\in GL_{in}(X)$ and note that if $\lambda \neq 1$, then
$T-\lambda \Id=(1-\lambda)(\Id+\frac{U}{1-\lambda})$ is Fredholm with index $0$ and so $T-\lambda \Id$ is a perturbation of an invertible operator by an operator in ${\mathcal F}(X)$. Therefore,
$\alpha(T-\lambda \Id)=\alpha(T)-\lambda \alpha(\Id)$ is an invertible element of ${\mathcal B}$ and hence $\lambda\notin \sigma\big(\alpha(T)\big)$.
We deduce that $\sigma\big(\alpha(T)\big)=\{1\}$, and since $\{\alpha(T)^n \del n \in \Z\}$ is bounded in the unital Banach algebra ${\mathcal B}$, Gelfand's theorem implies that
$\alpha(T)=\alpha(\Id)$. So $T-\Id$ belongs to $\AF(X)$, which concludes the proof of the complex case. Note that our proof in fact applies to any ideal ${\mathcal U}$ containing the  finite-rank operators and  such that any ${\mathcal U}$-perturbation of $\Id$ is Fredholm.

If instead  $X$ is real, we consider its complexification $\hat{X}$ and the ideal 
$$
{\mathcal U}=\In(X)+i\In(X)$$ of ${\mathcal L}(\hat{X})$, and observe that it contains ${\mathcal F}(\hat{X})={\mathcal F}(X)+i{\mathcal F}(X)$. 

We  claim that for all $U+iV$ in ${\mathcal U}$, $\Id+U+iV$ is Fredholm on $\hat{X}$.
Admitting the claim, we see that given $T \in GL_{in}(X)$, one can apply the proof in the complex case to $\hat{T}$, which is a ${\mathcal U}$-perturbation of $\Id$, and since then $\hat{T} \in {GL_{af}(\hat{X})}$ deduce that $T \in {GL_{af}(X)}$, thereby concluding the proof.

To prove the claim, note that since $\Id+U$ is Fredholm with index $0$ on $X$,  there exist $A \in GL(X)$ and $F \in {\mathcal F}(X)$ such that $\Id+U=A+F$. Then
$$\Id+U+iV=A+F+iV=A(\Id+A^{-1}F+iA^{-1}V),$$
which indicates that it is enough to prove that $\Id+iV$ is Fredholm for any $V \in {\mathcal U}$.

Fix such a $V$, write $\Id+V^2=B+L$, where $B \in GL(X)$ and $L \in {\mathcal F}(X)$, let
$F$ be the finite-dimensional subspace $B^{-1}LX$, let $H$ be a closed subspace such that
$X=F \oplus H$, let $\delta=d({S}_{H+iH},F+iF)>0$, and let $\epsilon>0$ be such that
$\sqrt{2}\epsilon\norm{B^{-1}}(1+\|V\|) <\delta$.

Let $x,y \in X$ and assume $\|(\Id+iV)(x+iy)\| \leqslant\epsilon$. An easy computation shows that
$$
\norm{x-Vy} \leqslant \epsilon
$$
and
$$
\norm{Vx+y} \leqslant \epsilon,
$$
whereby 
$$
\norm{Bx+Lx}=\norm{(\Id+V^2)x} \leqslant \epsilon(1+\norm{V}),
$$
and
$$
d(x,F) \leqslant \epsilon \norm{B^{-1}}(1+\norm{V}) < \delta/\sqrt{2}.
$$
Similarly $d(y,F) < \delta/\sqrt{2}$, and  so $d(x+iy,F+iF) < \delta$.
Conversely, this means that if $x+iy$ is a norm one vector in $H+iH$, then 
$$
\|(\Id+iV)(x+iy)\| > \epsilon,
$$
and so the restriction of $\Id+iV$ to the finite codimensional subspace $H+iH$ is an isomorphism onto its image. This proves that $\Id+iV$ has finite dimensional kernel and closed image. In particular, its Fredholm index is defined, with the possible value $-\infty$, but then the continuity of the index implies that this index is $0$ and therefore that $\Id+iV$ is Fredholm. This concludes the proof of the claim and of the theorem.
\end{proof}

Note that if $\K$ denotes the  scalar field of $X$, then 
$$
\K^\star=\{\lambda \Id\del \lambda\in \K\setminus \{0\}\}
$$
is a norm-closed subgroup of $GL(X)$. Also, if $\ku I$ is a  proper ideal in $\ku L(X)$, then $\K^\star\cap GL_{\ku I}(X)=\{\Id\}$ and so the group 
$$
\{\lambda \Id+A\in GL(X)\del \lambda\in \K\setminus \{0\}\text{ and }A\in \ku I\;\}
$$
of non-zero scalar multiples of elements of $GL_{\ku I}(X)$ splits as a direct product 
$$
\K^\star\times GL_{\ku I}(X).
$$
Moreover, since then $\ku J=\overline{\ku I}$ is a norm closed proper ideal, both $\K^\star$ and $GL_{\ku J}(X)$ are norm closed  and so the decompositions $\K^\star\times GL_{\ku J}(X)$ and hence also $\K^\star\times GL_{\ku I}(X)$ are topological direct products with respect to the norm  topology. In particular, this applies to the ideals $\ku F$ and $\ku {A\!F}$.
So though our ultimate interest is in, e.g.,  the group $\K^\star \times GL_f(X)$, in many situations this splitting allows us to focus on only the non-trivial part, namely $GL_f(X)$.

\begin{prop}\label{topos}
Suppose $X$ is a Banach space with separable dual and  $\pi\colon G\til  \K^\star  \times GL_{in}(X)$ is a bounded, weakly continuous representation of a Polish group $G$. Then $\pi$ is norm continuous.

It follows that if $G \leqslant \K^\star  \times GL_{in}(X)$   is a bounded subgroup,  strongly closed in $GL(X)$,
then the strong operator and  norm topologies coincide on $G$.
\end{prop}

\begin{proof} Composing $\pi$ with the coordinate projection from $\K^\star \times  GL_{in}(X)$ onto $GL_{in}(X)$ and using Theorem \ref{in then af}, we see that the representation has image in $\K^\star  \times GL_{af}(X)$.
We remark that, since $X^*$ is separable, the ideal $\ku F(X)$ of finite-rank operators on $X$ is separable for the norm topology, whence also $\AF(X)=\overline {\ku F(X)}^{\norm\cdot}$ and $\K^\star \times  GL_{af}(X)$  are norm separable. Moreover, by Proposition \ref{weakstrong}, $\pi$ is strongly continuous, whereby, for every $\eps>0$ and $x\in X$, the set 
$$
\ku U_{\eps,x}=\{g\in G\del \|\pi(g)x-x\|<\eps\}
$$
is open in $G$. Thus, if $\{x_n\}_{n\in \N}\subseteq X$ is dense in the unit ball of $X$, we see that
$$
\{g\in G\del \|\pi(g)-\Id\|<\eps\}=\bigcup_{m\geqslant 1}\bigcap_{n\in \N}\ku U_{\eps-1/m, x_n}
$$
is Borel in $G$. So $\pi$ is a Borel measurable homomorphism from a Polish group to a norm separable topological group and therefore is norm continuous by Pettis' Theorem (see \cite{kechris}, (9.10)).

If now instead $G\leqslant \K^\star  \times GL_{in}(X)$ is a bounded subgroup, strongly closed in $GL(X)$, then $G$ is Polish in the strong operator topology and so the tautological representation on $X$ is norm continuous, implying that every norm open set in $G$ is also strongly open. It follows that the two topologies coincide on $G$.
\end{proof}

Observe that if a space $X$ has an unconditional basis and $G$ is the bounded group of isomorphisms acting by change of signs of the coordinates on the basis, then $G$ is an uncountable discrete group in the norm topology and is just the Cantor group in the strong operator topology.
So there is no hope of extending  Proposition \ref{topos} to arbitrary strongly closed bounded subgroups $G\leqslant GL(X)$ when $X$ has an unconditional basis, and, in many cases, to ensure the norm-separability of any bounded $G\leqslant GL(X)$, we shall even have to assume that $X$ does not contain any unconditional basic sequences.

\subsection{Near triviality}As a corollary of Proposition \ref{topos} and Theorem \ref{NSS}, one sees that if $X$ is a Banach space with separable dual and  $G\leqslant {GL_{in}(X)}$ is compact in the strong operator topology, then $G$ is a compact Lie group. However,  we can prove an even stronger result  that also allows us to bypass the result of Gleason--Montgomery--Yamabe--Zippin. The central notion here is that of near triviality.

\begin{defi}
Let $X$ be a Banach space and $G\leqslant GL(X)$ a subgroup. We say that $G$ acts {\em nearly trivially} on $X$ if $X$ admits a decomposition into $G$-invariant subspaces,
$$
X=H\oplus F,
$$
such that $F$ is finite-dimensional and for all $T \in G$ there exists $\lambda_T$ such that $T|_H=\lambda_T \Id_H$. 
\end{defi}

We remark that a subgroup $G\leqslant GL(X)$ acts nearly trivially on $X$ if and only if the subgroup $c(G)$ of $GL(\hat{X})$ acts nearly trivially on the complexification $\hat{X}$.

Note that, when $G$ acts nearly trivially on $X$, the strong operator topology on $G$ is just the topology of pointwise convergence on $F'=F \oplus L$, where $L$ is an arbitrary one dimensional subspace of $H$ (or $L=\{0\}$ if $H$ is trivial), and so $T\in G\mapsto T|_{F'}\in GL(F')$ is a topological group embedding. Since any strongly closed bounded subgroup of $GL(F')$ is a compact Lie group, it follows that if $G$ is strongly closed and bounded in $GL(X)$, then $G$ is also a compact Lie group.

We also remark that in this case one has $\im(T-\lambda_T\Id)\subseteq F$ for all $T\in G$. But, in fact, this observation leads to the following equivalent characterisation of near triviality for bounded subgroups, which, for simplicity, we only state for subgroups of $GL_f(X)$.

\begin{lemme}\label{near triviality}
Let $X$ be a Banach space and $G\leqslant GL_f(X)$ a bounded subgroup. Then the following are equivalent.
\begin{enumerate}
\item $G$ acts nearly trivially on $X$,
\item there is a finite-dimensional $F\subseteq X$ such that $\im(T-\Id)\subseteq F$ for all $T\in G$.
\end{enumerate}

Moreover, in this case, 
$$
X=\bigcap_{T\in G}\ker (T-\Id)\oplus {\rm span}\big(\bigcup_{T\in G}\im (T-\Id)\big)
$$
is a decomposition witnessing near triviality.
\end{lemme}

\begin{proof}
One direction has already been noted, so suppose instead that (2) holds and let $F\subseteq X$ be the finite-dimensional subspace $F={\rm span}\big(\bigcup_{T\in G}\im (T-\Id)\big)$. Then for all $x\in X$ and $T\in G$, $Tx=x+f$ for some $f\in F$ with $\|f\|\leqslant \|T-\Id\|\cdot \norm x$. So, for any $x\in X$, $G\cdot x\subseteq x+(2\norm G\norm x) B_F$, showing that the orbit of $x$ is relatively compact. Moreover, $F$ is $G$-invariant and if $Y\subseteq X$ is any $G$-invariant subspace with $Y\cap F=\{0\}$, then $Tx=x$ for all $x\in Y$ and $T\in G$.

Since $G$ is almost periodic, by Proposition \ref{Gcompact}, $X$ is the closed linear space of its finite-dimensional irreducible subspaces. Therefore, as $Y\subseteq \bigcap_{T\in G}\ker (T-\Id)$ for any irreducible $Y\subseteq X$ with $Y\cap F=\{0\}$, we see that $X=F\oplus \bigcap_{T\in G}\ker (T-\Id)$, which finishes the proof.
\end{proof}

As easy applications, we have the following lemmas.

\begin{lemme}\label{finitely generated}
Suppose $X$ is a Banach space and $G\leqslant GL_f(X)$ is a finitely generated bounded subgroup. Then $G$ acts nearly trivially on $X$.  
\end{lemme}

\begin{proof}
Let $G=\langle T_1,\ldots,T_n\rangle$ and put $F=\im(T_1-\Id)+\ldots+\im(T_n-\Id)$, which is finite-dimensional. 
Note now that $T_i\inv-\Id=(T_i-\Id)(-T_i\inv)$, and so $\im(T_i\inv-\Id)\subseteq \im(T_i-\Id)\subseteq F$. Moreover, for $T,S\in GL(X)$, 
$$
TS-\Id=(S-\Id)+(T-\Id)S,
$$
and so, if $\im (T-\Id)\subseteq F$ and $\im(S-\Id)\subseteq F$, then also $\im(TS-\Id)\subseteq F$. It thus follows that $\im(T-\Id)\subseteq F$ for all $T\in G$, whence Lemma \ref{near triviality} applies.
\end{proof}

\begin{lemme}\label{single generator}
Suppose $X$ is a Banach space and $T\in GL_f(X)$ is an isometry. Then
$$
X=\ker (T-\Id)\oplus \im(T-\Id).
$$
Moreover, if $X$ is complex,  there are eigenvectors $x_i$ such that $ \im(T-\Id)=[x_1,\ldots,x_n]$.
\end{lemme}

\begin{proof}
As in Lemma \ref{finitely generated}, we see that $\im(T^n-\Id)\subseteq \im(T-\Id)$ for all $n\in \Z$. The result now follows from Lemma \ref{near triviality}. The moreover part follows from the fact that any isometry of a finite-dimensional complex space can be diagonalised.
\end{proof}

\begin{thm}\label{christian}
Suppose $X$ is a Banach space with separable dual and $G\leqslant GL_{in}(X)$ is an almost periodic subgroup, strongly closed in $GL(X)$. Then $G$ acts nearly trivially on $X$ and hence is a compact Lie group.
\end{thm}

\begin{proof}
There are several ways of proving this, e.g., one based on the structure theory for norm continuous representations (cf.  \cite{shtern}) of compact groups. But we will give a simple direct argument as follows.

Since $G$ is almost periodic, by Proposition \ref{Gcompact}, $G$ is compact in the strong operator topology and $X$ is the closed linear span of its finite-dimensional $G$-invariant subspaces. So, as $X$ is separable,  by taking finite sums of these we can find an increasing sequence $F_0\subseteq F_1\subseteq F_2\subseteq \ldots \subseteq X$ of finite-dimensional $G$-invariant subspaces such that $X=\overline{\bigcup_{n\geqslant 0}F_n}$. Moreover, by Proposition \ref{topos}, the norm and strong operator topologies coincide on $G$.

Let $\go A\subseteq \ku L(X)$ be the subalgebra generated by $G$ and define, for every $n$, the unital algebra homomorphism
$\pi_n\colon \go A\til \ku L(X/F_n)$ by 
$$
\pi_n(A)\big(x+F_n\big)=Ax+F_n.
$$
Note that $\norm{\pi_n(A)}\leqslant \norm{\pi_m(A)}\leqslant \norm A$ for all $n\geqslant m$.

We claim that there is an $n$ such that $\pi_n(T-\Id)=0$ for all $T\in G$. To see this,  assume the contrary and note that, by Theorem \ref{NSS}, for every $n$ there is some $T_n\in G$ such that $\norm{\pi_n(T_n-\Id)}=\norm{\pi_n(T_n) -\Id}\geqslant \sqrt 2$. So for $m\leqslant n$ we have 
$$
\norm{\pi_m(T_n-\Id)}\geqslant \norm{\pi_n(T_n-\Id)}\geqslant \sqrt 2.
$$
Moreover, by passing to a subsequence, we can suppose that the $T_n$ converge in norm to some $T\in G$, whence $\norm{\pi_m(T-\Id)}\geqslant \sqrt 2$ for all $m$. 

Now, by Theorem \ref{in then af}, $G\leqslant GL_{af}(X)$, so $A=T-\Id\in \ku{A\!F}$. Since $A$ is a norm limit of finite-rank operators, there is a finite-dimensional subspace $F\subseteq X$ such that $\norm{Ax+F}<1$ for all $x\in X$, $\norm x\leqslant 1$. Using that the $F_n$ are increasing and $X=\overline{\bigcup_{n\geqslant 0}F_n}$, we see that there is an $n$ such that
$$
\sqrt 2\leqslant \norm{\pi_n(A)}=\sup_{\norm{x}=1}\norm{Ax+F_n}\leqslant\sup_{\norm{x}=1}\norm{Ax+F}+1/3\leqslant 4/3,
$$
which is absurd. So fix an $n$ such that  $\pi_n(T-\Id)=0$  and thus $\im (T-\Id)\subseteq F_n$ for all $T\in G$. The result now follows from Lemma \ref{near triviality}.
\end{proof}

To simplify notation, we let
$$
{\rm Isom}_f(X)={\rm Isom}(X)\cap GL_f(X)
$$ 
and
$$
{\rm Isom}_{af}(X)={\rm Isom}(X)\cap GL_{af}(X)
$$ 
denote the normal subgroups of ${\rm Isom}(X)$ consisting of all so-called {\em finite-dimensional}, resp. {\em almost finite-dimensional} isometries. Recall that by  Theorem \ref{in then af}, any isometry which is an inessential perturbation of $\Id$ must belong to
${\rm Isom}_{af}(X)$. The next lemma shows that under additional conditions one may replace
${\rm Isom}_{af}(X)$ by ${\rm Isom}_f(X)$.

\begin{lemme}\label{gelfand}
Let $X$ be a complex Banach space and let $T \in {\rm Isom}_{af}(X)$. If $T$ has finite spectrum, then $T \in {\rm Isom}_f(X)$.
\end{lemme}

\begin{proof} Write $\sigma(T)=\{1,\lambda_1,\ldots,\lambda_n\}$ and let $P$ be the spectral projection of $T$ corresponding to the spectral set $\{1\}$.
Then $T(PX) = PX$, $\sigma(T|_{PX})=\{1\}$ and $\sup_{n\in \Z}\norm{(T|_{PX})^n}<\infty$, so, by the result of Gelfand, $T|_{PX}=\Id$. By the same reasoning $T|_{X_i}=\lambda_i \Id_{X_i}$, where $X_i$ is the range of the spectral projection associated to $\{\lambda_i\}$ for $i=1,\ldots,n$.
Since $T$ is an almost finite-rank perturbation of the identity, all elements of $\sigma(T)$ different from $1$ have finite multiplicity and therefore $PX$ has finite-codimension.
\end{proof}

As in the case of $GL_f(X)$, it will  often be enough to study the group ${\rm Isom}_f(X)$, although our interest 
will really be in the subgroup $\{-1,1\}\times {\rm Isom}_f(X)$ of ${\rm Isom}(X)$ (respectively
$\T \times {\rm Isom}_f(X)$ in the complex case). The same holds in relation to
${\rm Isom}_{af}(X)$ and $\{-1,1\}\times{\rm Isom}_{af}(X)$
(respectively $\T \times {\rm Isom}_{af}(X)$).


\section{Decompositions of separable reflexive spaces by isometries}

\subsection{Duality mappings}
Let $X$ be a Banach space. Recall that a {\em support functional} for $x \in   S_X$  is a  functional $\phi$ in $  S_{X^*}$ such that $\phi(x)=1$. Support functionals always exist by the Hahn--Banach theorem.
We shall denote by $Jx$ the set of support functionals of $x \in   S_X$ and extend  $J$ to all of $X$ by positive homogeneity, that is, $J(tx)=tJx$ for all $t \geqslant 0$ and $x \in   S_X$. Also, for $Y\subseteq X$, $J[Y]$ denotes the set of support functionals for $x\in Y$.

\begin{lemme}\label{preorthogonal}
Let $X$ be a Banach space and $Y$  a closed linear subspace of $X$. Then the following hold.
\begin{itemize}
\item[(a)] $Y$ and $J[Y]_\perp$ form a direct sum in $X$ and the corresponding projection from the subspace
$$
Y \oplus J[Y]_{\perp}
$$
onto $Y$ has norm at most $1$.
\item[(b)] If $Y$ is  reflexive and   $J[Y]$ is a closed linear subspace of $X^*$, then
$$
X^*= J[Y]\oplus Y^{\perp}
$$
and the corresponding projection from $X^*$ onto $J[Y]$ has norm at most $1$.
\end{itemize}
\end{lemme}

\begin{proof}
Suppose $y\in   S_Y$, $z\in J[Y]_{\perp}$, and let $\phi\in Jy$. Then
$$
\norm{y+z} \geqslant \phi(y+z)=\phi(y)=1=\norm{y},
$$
which implies (a). A similar argument shows that if $J[Y]$ is a closed linear subspace of $X^*$, then $J[Y]$ and $Y^{\perp}$ form a direct sum in $X^*$ and the corresponding projection from the subspace $$J[Y] \oplus Y^{\perp}$$ onto $J[Y]$ has norm at most $1$. 

For (b), assume that $Y$ is  reflexive  and $J[Y]$ is a closed linear subspace of $X^*$. To see that $X^*=J[Y] \oplus Y^{\perp}$, fix $\psi\in X^*$ and let $\phi \in X^*$ be a Hahn--Banach extension of $\psi|_Y$ to all of $X$ with $\|\phi\|=\|\psi|_Y\|$. Then $\psi|_Y=\phi|_Y$, whence $\psi-\phi\in Y^\perp$ and  $\|\phi\|=\|\psi|_Y\|=\|\phi|_Y\|$. On the other hand, since $\|\phi\|=\|\phi|_Y\|$ and $Y$ is reflexive, $\phi|_Y$ and thus $\phi$ attain their norm on $Y$, which means that $\phi=Jy$ for some $y\in Y$, i.e., $\phi\in J[Y]$. So $\psi=\phi+(\psi-\phi)\in J[Y]\oplus Y^\perp$.
\end{proof}

We refer the reader to \cite{DGZ} for more general results in this direction, see, for example, Lemma 2.4 p. 239 for information about the class of Weakly Countably Determined spaces.

Recall that a norm $\norm\cdot$  on a Banach space $X$ is {\em  G\^ateaux differentiable} if for every $x \in   S_X$ and $h \in X$,
$$
\lim_{t \rightarrow 0,\; t \in \R}\frac{\|x+th\|-\|x\|}{t}
$$
exists and is a linear continuous function in $h$. We note that this only depends on the $\R$-linear structure of $X$.
When the norm on $X$ is G\^ateaux differentiable, the support functional is unique for all $x \in   S_X$. This is proved in \cite{DGZ} in the real case and is also true in the complex case, as observed in Section \ref{complexification}. So, provided that the norm is G\^ateaux differentiable,  $Jx$ is a singleton for all $x \in   S_X$ and we can therefore  see $J$ as a map from $X$ into $X^*$.

Therefore, assuming that the norms on $X$ and $X^*$ are both G\^ateaux differentiable, the duality map is defined from $X$ to $X^*$ and from $X^*$ to $X^{**}$, where to avoid confusion we denote the second by $J_*\colon X^*\til X^{**}$. 
The following two lemmas are now almost immediate from the definition of $J$ and $J_*$.

\begin{lemme}\label{dualityJ}
Let $X$ be a reflexive space with a G\^ateaux differentiable norm, whose dual norm is G\^ateaux differentiable. Then $J\colon X\til X^*$ is a bijection with inverse $J_*\colon X^*\til X$.
\end{lemme}

\begin{lemme}\label{facto}
Let $X$ have a G\^ateaux differentiable norm and let $T$ be an isometry of $X$. Then
the maps $J T\inv $ and $T^{*}J$ coincide on $X$.
\end{lemme}

\begin{proof}
For any $x \in   S_X$, $T^{*}(Jx)(T\inv x)=1$, which shows that
$J(T\inv x)=T^{*}(Jx)$. This extends to all of $X$ by positive homogeneity.
\end{proof}


\subsection{LUR renormings and isometries}
Recall that a norm on a Banach space $X$ is said to be {\em locally uniformly rotund} or, in short, {\em LUR}, if
$$
\a x_0 \in   S_X\;\a \epsilon>0\; \e \delta>0\; \a x \in   S_X
\;\big( \norm{x-x_0} \geqslant \epsilon\saa\norm{x+x_0} \leqslant 2-\delta\big).
$$
For other characterizations of LUR norms we refer to \cite{DGZ}. Any LUR norm is strictly convex, meaning that the associated closed unit ball is strictly convex. Note that this definition does not depend on $X$ being seen as real or complex.

\begin{prop}\label{maps}
Assume $X$ is a Banach space whose dual norm is LUR. Then the duality mapping $J\colon X\til X^*$
is well-defined and norm-continuous.

Therefore, if $X$ is reflexive and both the norm and the dual norm are LUR, then $J\colon X\til X^*$ is a homeomorphism with inverse $J_*\colon X^*\til X$.
\end{prop}

\begin{proof} It follows from the LUR property in $X^*$ that the norm on $X$ is G\^ateaux-differentiable \cite{DGZ} and therefore that $J$ is well-defined. Now given $x_0 \in   S_X$ and $\epsilon>0$, let $\delta>0$ be associated to $\epsilon$ by the LUR-property at $\phi_0=Jx_0$.
If for some $x \in   S_X$, $\|Jx-Jx_0\| \geqslant \epsilon$, then
$\|Jx+Jx_0\| \leqslant 2-\delta$ and therefore
$$
|1-(Jx)(x_0)|=|2-(Jx+Jx_0)(x_0)| \geqslant 2-\|Jx+Jx_0\| \geqslant \delta.
$$
Thus
$$
\norm{x-x_0} \geqslant |(Jx)(x-x_0)|=|1-(Jx)(x_0)| \geqslant \delta.
$$
This proves that $J$ is continuous as a map from $  S_X$ to $  S_{X^*}$ and therefore from  $X$ to $X^*$.
\end{proof}

It is a well-known result of renorming theory, due to M.I. Kadec, that any separable real space admits an equivalent LUR norm, see \cite{DGZ} Chapter II. However, we are interested in LUR renormings which, in some sense, keep track of the original group of isometries on the space. The next proposition is the first of a series of results in that direction. 

For expositional ease, if $\norm{\cdot}$ is a norm on a Banach space $X$, we shall denote the induced norm on the dual space $X^*$ by $\norm{\cdot}^*$.

\begin{prop}\label{renormagebaire} 
Let $X$ be a Banach space and $G\leqslant GL(X)$ a bounded subgroup. Assume that $X$ admits $G$-invariant equivalent norms $\norm{\cdot}_0$ and $\norm{\cdot}_1$ such that $\norm{\cdot}_0$ and $\norm{\cdot}_1^*$ are LUR. Then $X$ admits a $G$-invariant equivalent norm $\norm{\cdot}_2$ such that both $\norm{\cdot}_2$ and $\norm{\cdot}_2^*$ are LUR.

Moreover, for any $G$-invariant equivalent norm $\norm{\cdot}$ on $X$ and $\eps>0$, one can choose $\norm{\cdot}_2$ to be $(1+\epsilon)$-equivalent to $\norm{\cdot}$. 
\end{prop}

\begin{proof} 
Let $\ku N$ denote the space of $G$-invariant equivalent norms on $X$ equipped with the complete metric 
$$
d\big(\norm{\cdot},\triple{\cdot}\big)=\sup_{x\neq 0}\Big|\log\frac {\norm{x}}{\triple{x}}\Big|.
$$ 
Let also $\ku L\subseteq \ku N$ denote the subset of all norms that are LUR. We first show that $\ku L$ is comeagre in $\ku N$.

Note first that $\norm{\cdot}_0\in \ku L$ and define for all $k\geqslant 1$ the open set 
$$
\ku L_k=\Big\{\triple{\cdot}\in \ku N\Del \e \norm{\cdot}\in \ku N\; \e n\geqslant k\;\;d\big(\triple{\cdot},\big(\norm{\cdot}^2+\frac 1n\norm{\cdot}_0^2\big)^{1/2} \big)<\frac 1{n^2}\Big\}
$$
Now, if $\norm{\cdot}\in \ku N$, then for all  $n\geqslant k$, $\big(\norm{\cdot}^2+\frac 1n\norm{\cdot}_0^2\big)^{1/2} \in \ku L_k$, while on 
the other hand $\big(\norm{\cdot}^2+\frac 1n\norm{\cdot}_0^2\big)^{1/2} \Lim{n\til \infty}\norm{\cdot}$, which shows that $\norm{\cdot}\in \overline{\ku L_k}$ for every $k\geqslant 1$. Thus, $\ku L_k$ is dense open for every $k$ and hence $\bigcap_{k\geqslant 1}\ku L_k$ is comeagre in $\ku N$ and, as shown in \cite{DGZ}, is a subset of $\ku L$. So $\ku L$ is comeagre in $\ku N$. 

Similarly, one shows, using that $\norm{\cdot}_1^*$ is LUR, that the set 
$$
\ku M=\big\{\norm{\cdot}\in \ku N\del \norm{\cdot}^* \text{ is LUR }\big\}
$$
is comeagre in $\ku N$. It thus suffices to choose $\norm{\cdot}_2$ in the comeagre and thus dense intersection $\ku L\cap \ku M$.
\end{proof}

The following result was proved by G. Lancien, see \cite{L} Theorem 2.1, Remark 1, p. 639,  Theorem 2.3,  and the observation and remark p. 640.

\begin{thm}[G. Lancien]\label{lancien}
Let $(X,\norm\cdot)$ be a separable Banach space and set $G={\rm Isom}(X,\norm{\cdot})$. Then the following hold.
\begin{itemize}
\item[(a)] If $X$ has the Radon--Nikodym Property, then $X$ admits an equivalent $G$-invariant LUR norm.
\item[(b)] If $X^*$ is separable,  then $X$ admits an equivalent  $G$-invariant norm whose dual norm is LUR.
\end{itemize}
\end{thm}

We should note here that the result of Lancien is proved and stated for real spaces, but also holds for complex spaces. Indeed, suppose $X$ is a complex space and let $X_\R$ denote the space $X$ seen as a Banach space over the real field.
Any equivalent real norm on $X_\R$, which is ${\rm Isom}(X_\R,\norm{\cdot})$-invariant, must be invariant under all isometries of the form $\lambda \Id$ for $\lambda \in \T$ and hence is actually a complex norm that is ${\rm Isom}(X,\norm{\cdot})$-invariant.
Thus, in order to obtain Lancien's result for a complex space $X$, and modulo the observations in Section \ref{complexification}, it suffices to simply apply it to $X_\R$.

Combining Proposition \ref{renormagebaire} and Theorem \ref{lancien}, plus the fact that any bounded subgroup of $GL(X)$ is an isometry group in some equivalent norm, we obtain the following result.

\begin{thm}\label{newlancien}
Let $X$ be a  Banach space with the Radon--Nikodym Property and separable dual and let $G\leqslant GL(X)$ be a bounded subgroup. Then $X$ admits an equivalent $G$-invariant LUR norm whose dual norm is also LUR. 
\end{thm}

This should be compared with the relatively immediate fact (see, e.g., Proposition 2.3 \cite{furman}) that any super-reflexive space admits an isometry invariant uniformly convex norm.

Since any reflexive space has the RNP \cite{BL},  the conclusion of Theorem \ref{newlancien} holds, in particular, for any separable reflexive space. As a first consequence of this, we have the following lemma.

\begin{lemme}\label{AB} 
Let $X$ be separable reflexive and $G$ be a bounded subgroup of $GL(X)$. Then for each $x \in X$, $\overline{\rm conv}(G\cdot x)$ contains a unique $G$-fixed point.
\end{lemme}

\begin{proof} By renorming using Theorem \ref{newlancien}, we may assume that $G$ is a group of isometries and that the norms on $X$ and $X^*$ are LUR. Fix $x \in X$ and let $C=\overline{\rm conv}(G\cdot x)$. By reflexivity, $C$ is  weakly compact and thus contains a point $x_0$ of minimal norm. Furthermore, by strict convexity, $x_0$ is unique and therefore fixed by $G$.

Assume $x_1$ is another $G$-fixed point in $C$ and let $\phi \in X^*$. As above, we may find
a unique functional $\phi_0$ of minimal norm in $\overline{\rm conv}(G\cdot \phi)$, which is therefore fixed by $G$. It follows that $\phi_0$ is constant on $G\cdot x$ and thus also on $C=\overline{\rm conv}(G\cdot x)$, whence $\phi_0(x_1)=\phi_0(x_0)=\phi_0(x)$ and  
$\phi_0(x_1-x_0)=0$. Now since $x_1-x_0$ is $G$-fixed, seen as a functional on $X^*$, it must be constant on $\overline{\rm conv}(G\cdot\phi)$ and so $\phi(x_1-x_0)=0$.
As $\phi$ was arbitrary, this proves that $x_1=x_0$.
\end{proof}


\subsection{Decompositions}
Suppose $G\leqslant GL(X)$ is a bounded group of automorphisms of a Banach space $X$ and define the following set of closed $G$-invariant subspaces
\begin{itemize}
\item $H_G=\{x\in X\del x\text{ is fixed by } G\}$,
\item $H_{G^*}=\{\phi\in X^*\del \phi\text{ is fixed by } G\}$,
\item $K_G=\{x\in X\del  x\text{ is almost periodic}\,\}$,
\item $K_{G^*}=\{\phi\in X^*\del \phi\text{ is almost periodic}\,\}$.
\end{itemize}
Note also that $H_G \subseteq K_G$ and $H_{G^*} \subseteq K_{G^*}$.
When $G$ is generated by a single element, i.e., $G=\langle T\rangle$, we shall simply write $H_T$, $H_{T^*}$, etc., instead of the more cumbersome $H_{\langle T\rangle}$ and $H_{\langle T\rangle^*}$. Define also the subspaces
\begin{itemize}
\item $F_T=\im (T-\Id)$,
\item $F_{T^*}=\im (T^*-\Id)$.
\end{itemize}
Then, as is easy to verify, $H_{T^*}=(F_T)^\perp$, $H_T=(F_{T^*})_\perp$ and $\overline{F_T}=(H_{T^*})_\perp$.

K. Jacobs \cite{jacobs} and later K. de Leeuw and I. Glicksberg \cite{deleeuw} studied the conditions under which the the subspace $K_G$ of almost periodic vectors in a reflexive space $X$ admits a $G$-invariant complement and were able to show this, e.g.,  in the case $X$ and $X^*$ are both strictly convex and $G$ is a group of isometries (see \cite{deleeuw} Corollary 4.14). Since, by Theorem \ref{newlancien}, when $X$ is separable reflexive and $G\leqslant GL(X)$ is a bounded subgroup, we can renorm $X$ so that these conditions hold, we see that one can apply their results in this setting. However, we can give a direct proof of this decomposition again using  Theorem \ref{newlancien} and the properties of the duality mapping, and also obtain quantitative estimates for the norm of the associated projections, which shall be needed later on.

By \cite{jacobs}, if $X$ is any Banach space, $G\leqslant GL(X)$ is a bounded subgroup and $x\in X$, we say that $x$ is  {\em furtive} if $0\in \overline{G\cdot x}^{\rm w}$. Of course, when $X^*$ is separable, this is equivalent to requiring that there is a sequence $T_n\in G$ such that $T_nx\Lim{\rm w}0$. Note that furtive vectors do not in general form a linear subspace of $X$.

\begin{lemme} Let $X$ be a Banach space and $G \leqslant GL(X)$ be a bounded subgroup. Then the following implications hold:
\[
\begin{tabular}{ccc}
$x$ is furtive & $\Rightarrow$ &  $x \in (K_{G^*})_\perp$ \\

$\Downarrow$&       &  $\Downarrow$\\ 

$0\in \overline{\rm conv}(G\cdot x)$ & $\Rightarrow$ & $x \in (H_{G^*})_\perp$  \\
\end{tabular}
\]
\end{lemme}

\begin{proof}
If $x\in X$ is furtive and $\phi\in K_{G^*}$, we claim that $\phi(x)=0$. To see this, assume without loss of generality that $G$ is a group of isometries, fix $\eps>0$ and pick an $\eps$-dense subset $\psi_1,\ldots, \psi_k$ of $G\cdot \phi$. Then, since $x$ is furtive, there is some $T\in G$ such that $|T^*\psi_i(x)|=|\psi_i(Tx)|<\eps$ for every $i=1,\ldots, k$. Picking $i$ such that $\|\phi-T^*\psi_i\|=\|(T\inv)^*\phi-\psi_i\|<\eps$, we see that 
$|\phi(x)|\leqslant \|\phi-T^*\psi_i\|\norm x+|T^*\psi_i(x)|<\eps\|x\|+\eps$. So, as $\eps>0$ is arbitrary, it follows that $\phi(x)=0$. In other words, any furtive vector belongs to $(K_{G^*})_\perp$.

Note also that if $\phi\in H_{G^*}$, $T_1,\ldots, T_n\in G$ and $\lambda_i\geqslant 0$ are such that $\sum_{i=1}^n\lambda_i=1$, then, for any $x\in X$,
$$
\phi(\sum_{i=1}^n\lambda_iT_ix)=\sum_{i=1}^n\lambda_iT^*_i\phi(x)=\sum_{i=1}^n\lambda_i\phi(x)=\phi(x).
$$
Hence, if $0\in \overline{\rm conv}(G\cdot x)$ and $\phi\in H_{G^*}$, then $\phi(x)=0$, showing the second implication. Finally, the vertical implications are trivial.
\end{proof}

It turns out that if $X$ is separable reflexive, then the horizontal implications  reverse, which will allow us to identify  the Jacobs - de Leeuw - Glicksberg type decomposition with a decomposition provided by the duality mapping.
 
\begin{thm}\label{decomp}
Let $X$ be a separable reflexive space and suppose $G\leqslant GL(X)$ is a bounded subgroup. Then $X$ admits the following $G$-invariant decompositions.

\begin{itemize}
\item[(a)](Alaoglu - Birkhoff type decomposition)
$$
X=H_G\oplus (H_{G^*})_\perp,
$$
where, if $\ku S\subseteq G$ generates $G$, then
$$
(H_{G^*})_\perp=\overline{\rm span}\big(\bigcup_{T\in \ku S}F_T\big)=\{x\in X\del 0\in \overline{\rm conv}(G\cdot x)\}.
$$
Moreover, the projection $P\colon X\til H_G$ is given by 
$$
Px=\text{ the unique point in }H_G\cap \overline{\rm conv}(G\cdot x)
$$
\item[(b)](Jacobs - de Leeuw - Glicksberg type decomposition)
$$
X=K_G\oplus (K_{G^*})_\perp,
$$
where
$$
(K_{G^*})_\perp=\{x\in X\del x\textrm{ is furtive}\}=\{x\in X\del \e T_n\in G\;\; T_nx\Lim{w}0\}.
$$
\end{itemize}
Moreover, the projections onto each summand have norm bounded by $2\norm{G}^2$ (or $\norm{G}^2$ when the summand is $H_G$ or $K_G$), and either $G$ is almost periodic, i.e., $X=K_G$, or the subspace $(K_{G^*})_\perp$ is infinite-dimen\-sional. 

In particular, for any isometry $T$ of $X$,
$$
X=H_T\oplus \overline{F_T}.
$$
\end{thm}

\begin{proof}
Renorming $X$ by  $\triple x=\sup_{T \in G}\|Tx\|$, we can suppose $G$ is a group of isometries. Moreover, by a further renorming using Theorem \ref{newlancien}, we can suppose both the norm and its dual are LUR. Furthermore, using the quantitative estimate of Proposition \ref{renormagebaire}, and fixing $\epsilon>0$, we can ensure that the resulting norm on $X$ is $(1+\epsilon)\norm{G}$-equivalent to the original norm. 

Note then that, by Propositions \ref{maps} and \ref{facto}, $J\colon X\til X^*$ is a homeomorphism  satisfying $JT\inv =T^*J$, whence $J[H_G]\subseteq H_{G^*}$ and $J[K_G]\subseteq K_{G^*}$. Similarly, 
the inverse $J_*$ satisfies $J_*[H_{G^*}]\subseteq H_G$ and $J_*[K_{G^*}]\subseteq K_G$, whence $J[H_G]=H_{G^*}$ and $J[K_G]=K_{G^*}$. It follows from Lemma \ref{preorthogonal} (b) that 
$$
X=H_G\oplus (H_{G^*})_\perp=K_G\oplus (K_{G^*})_\perp,
$$
where the corresponding projections have norm at most $2$, or $1$ for the projections onto $H_G$ and $K_G$. Since $H_G\subseteq K_G$ and $H_{G^*}\subseteq 
K_{G^*}$, we see that the decompositions refine to $X=H_G\oplus\big((H_{G^*})_\perp\cap K_G\big)\oplus (K_{G^*})_\perp$. Since $\epsilon$ was arbitrary, the estimate on the norms of the projections in the original space follows immediately. 

Note that, since no non-zero $G$-orbit on  $(K_{G^*})_\perp$ is relatively compact, the latter space must either be infinite-dimensional or reduce to $\{0\}$.

Moreover, if $\ku S\subseteq G$ generates $G$, then 
$$
H_{G^*}=\bigcap_{T\in \ku S}H_{T^*}=\bigcap_{T\in \ku S}(F_T)^\perp=\big({\rm span}\big(\bigcup_{T\in \ku S}F_T\big)\big)^\perp
$$
and so $(H_{G^*})_\perp=\overline{\rm span}\big(\bigcup_{T\in \ku S}F_T\big)$.

Now, to see that $(K_{G^*})_\perp$ is the set of furtive vectors, note that since $G$ is a group of isometries of $(X,\triple\cdot)$ and both $\triple\cdot$ and $\triple{\cdot}^*$ are strictly convex, by Corollary 4.14 of \cite{deleeuw}, the set of furtive vectors form a linear subspace $Y$ such that $X=K_G\oplus Y$. Since, as we have seen, $Y\subseteq (K_{G^*})_\perp$ and also $X=K_G\oplus (K_{G^*})_\perp$, it follows that $Y=(K_{G^*})_\perp$.

Note now that for any $z\in (H_{G^*})_\perp$, we have $\overline{\rm conv}(G \cdot z)\subseteq (H_{G^*})_\perp$ and so the unique point in $H_G\cap \overline{\rm conv}(G \cdot z)$, that exists by Lemma \ref{AB}, must be $0$. Thus, if $y\in H_G$ and $z\in (H_{G^*})_\perp$, then for any $\eps>0$ there are $T_i\in G$ and $\lambda_i\geqslant 1$ with $\sum_{i=1}^n\lambda_i=1$ and 
$$
\norm{y-\sum_{i=1}^n\lambda_iT_i(y+z)}=\norm{y-\big(\sum_{i=1}^n\lambda_iy+\sum_{i=1}^n\lambda_iT_iz\big)}=\norm{\sum_{i=1}^n\lambda_iT_iz}<\eps.
$$
Since $\eps>0$ is arbitrary, this shows that the projection of $x=y+z$ onto $H_G$, namely $y$, belongs to $H_G\cap \overline{\rm conv}(G \cdot x)$, which, by Lemma \ref{AB}, implies that the projection $P$ of $X$ onto $H_G$ is given by $\{Px\}=H_G\cap \overline{\rm conv}(G\cdot x)$ and hence that $(H_{G^*})_\perp=\{x\in X\del 0\in \overline{\rm conv}(G\cdot x)\}$.
\end{proof}

We should mention that the Jacobs - de Leeuw - Glicksberg decomposition can be obtained for any weakly almost periodic subgroup $G\leqslant GL(X)$, where $X$ is separable. This can, e.g., be done by applying the Davis - Figiel - Johnson - Pe\l czy\'nski interpolation method \cite{davis} to Theorem \ref{decomp} to show that any vector in $X$ can be written as a sum of an almost periodic and a furtive vector. Moreover, by a result of S. Goldberg and P. Irwin \cite{goldberg}, the set of furtive vectors in $X$ is a closed linear subspace of $X$ forming a direct sum with $K_G$, and thus the result follows. For weaker results, but in the more general context of  a contractive semigroup $G$, one can consult \cite{krengel}.


\subsection{No totally bounded orbit}
By Theorem \ref{decomp}, the study of the isometry group of a separable, reflexive Banach space essentially reduces to the study of two separate cases, namely when all orbits are totally bounded and when no orbit is totally bounded. We shall now treat the second case.

\begin{lemme}\label{orbite non-compacte}
Suppose $G$ is a group acting by isometries on a complete metric space $M$. Then the following conditions are equivalent.
\begin{enumerate}
\item no orbit is totally bounded,
\item for every compact $C\subseteq M$, there is a $g\in G$ such that $g[C]\cap C=\tom$,
\item for any compact $C\subseteq M$ and $g\in G$ there are $f_1,f_2\in G$ such that $g=f_1f_2$ and $f_i[C]\cap C=\tom$.
\end{enumerate}
\end{lemme}

\begin{proof}
That (3) implies (1) is trivial. For the implication from (1) to (2), we shall use the following well-known fact from group theory.

\begin{claim}\label{neumann}
Suppose $F_1E_1\cup\ldots\cup F_nE_n$ is a covering of a group $G$, where  $E_i,F_i$ are subsets of $G$ with $F_i$ finite. Then there is a finite set $F\subseteq G$ and an $i$ such that $G=FE_iE_i\inv$.
\end{claim}

\begin{proof}
The proof is by induction on $n$, the case $n=1$ being trivial. So suppose the result holds for $n-1$ and let $G=F_1E_1\cup\ldots\cup F_nE_n$ be a covering.
Then, if $G\neq F_1E_1E_1\inv$, pick $g\in G\setminus F_1E_1E_1\inv$, whereby $gE_1\cap F_1E_1=\tom$. It follows that $gE_1\subseteq F_2E_2\cup\ldots\cup F_nE_n$ and hence
$$
F_1E_1=F_1g\inv \cdot gE_1\subseteq F_1g\inv F_2E_2\cup\ldots\cup F_1g\inv  F_nE_n.
$$
Thus,
$$
G=(F_1g\inv F_2\cup F_2)E_2\cup\ldots\cup (F_1g\inv  F_n\cup F_n)E_n,
$$
finishing the inductive step.
\end{proof}
Assume now that (1) holds and, for any $x\in M$ and $\delta>0$, define $V(x,\delta)=\{g\in G\del d(gx,x)<\delta\}$ and note that for any $f\in G$,
$$
f\cdot V(x,\delta)\cdot V(x,\delta)\inv=f\cdot V(x,\delta)\cdot V(x,\delta)\subseteq f\cdot V(x,2\delta)= \{g\in G\del d(gx,fx)<2\delta\}.
$$
Suppose $C\subseteq M$ is compact and assume towards a contradiction that $g[C]\cap C\neq\tom$ for any $g\in G$. Since no orbit is totally bounded, pick $\delta>0$ such that the orbit of no point of $C$ admits a finite covering by sets of diameter $8\delta$ and let
$$
B(x_1,\delta)\cup\ldots\cup B(x_n,\delta)
$$
be a finite subcover of the cover $\bigcup_{x\in C}B(x,\delta)$ of $C$. Then for any $g$ there are $i,j$ such that $g[B(x_i,\delta)]\cap B(x_j,\delta)\neq \tom$, whence $d(gx_i,x_j)<2\delta$.
Now, for every $i,j$ pick if possible some $f\in G$ such that $d(fx_i,x_j)<2\delta$ and let $F\subseteq G$ be the finite set of such $f$.  Then
$$
G= F\cdot V(x_1,4\delta)\cup \ldots \cup F\cdot V(x_n,4\delta),
$$
and so by Claim \ref{neumann} there is some finite $E\subseteq G$ and $i$ such that $G=E\cdot V(x_i,8\delta)$, whence $G\cdot x_i$ is covered by finitely many open balls of diameter $8\delta$, contradicting the choice of $\delta$.

Now, to see that (2) implies (3), suppose that $g\in G$ and $C\subseteq M$ is compact. Then there is some $f_1\in G$ such that
$$
f_1\big[C\cup g[C] \big]\cap \big(C\cup g[C]\big)=\tom,
$$
whence, in particular,  $f_1[C]\cap C=\tom$. Letting $f_2=f_1\inv g$, we have
\[\begin{split}
f_2[C]\cap C=f_1\inv g[C]\cap C\subseteq f_1\inv \big[C\cup g[C]\big]\cap \big(C\cup g[C]\big)=\tom.
\end{split}\]
Since $g=f_1f_2$, this finishes the proof.
\end{proof}

The following lemma is an immediate corollary.

\begin{lemme}\label{loin}
Let $X$ be a Banach space and $G$ a group of isometries of $X$ such that no non-zero $G$-orbit is totally bounded. Then for any finite-dimensional subspace $F$ of $X$ there exists a $T$ in $G$ such that $F \cap T[F]=\{0\}$. In fact, for any finite-dimensional subspace $F$ of $X$, any isometry $S$ in $G$ can be written as a product of two such $T$ in $G$.
\end{lemme}

With this in hand, we can now prove Theorem \ref{schauder}.

\begin{thm}\label{schauder}
Let $X$ be a  separable, reflexive  space and $G\leqslant GL(X)$ be a bounded subgroup such that no non-zero $G$-orbit is totally bounded. Assume also $G$ contains a non-trivial finite-rank perturbation of the identity. Then $X$ has a complemented subspace with a Schauder basis.
\end{thm}

\begin{proof}Note that for any $T,S\in G$, $F_{STS\inv}=S[F_T]$ and fix some $\Id\neq T_1\in G$ such that $F_{T_1}=\im(T_1-\Id)$ is a finite-dimensional subspace of $X$.

Since by renorming we may assume that $G$ is a group of isometries, Lemma \ref{loin} applies, and we can inductively define $S_2, S_3, \ldots\in G$ such that setting $T_n=S_nT_1S_n\inv$, we have 
$$
F_{T_n}\cap \big( F_{T_1}+F_{T_2}+\ldots+F_{T_{n-1}}\big)=S_n[F_{T_1}]\cap \big(F_{T_1}+S_2[F_{T_1}]+\ldots+S_{n-1}[F_{T_1}]\big)=\{0\}, 
$$
whence, in particular,  $\dim  \big( {F_{T_1}+\ldots+F_{T_{n}}}\big)=n\cdot \dim F_{T_1}$.
Moreover, by Theorem \ref{decomp}, we see that for every $n\geqslant 1$,
$$
X=\big(F_{T_1}+\ldots+F_{T_n}\big)\oplus H_{\langle T_1,\ldots,T_n\rangle},
$$
where the corresponding projections $\|P_n\|$ of $X$ onto the summands $\big(F_{T_1}+\ldots+F_{T_n}\big)$  are uniformly bounded by $2\|G\|^2$. Similarly, 
$$
X=\overline{\rm span}\big(\bigcup_{n\geqslant 1}F_{T_n}\big)\oplus H_{\langle T_1,T_2,\ldots\rangle}.
$$
Also, for $n\leqslant m$, we have $\big(F_{T_1}+\ldots+F_{T_n}\big)\subseteq \big(F_{T_1}+\ldots+F_{T_m}\big)$ and $H_{\langle T_1,\ldots,T_n\rangle}\supseteq H_{\langle T_1,\ldots,T_m\rangle}$, from which it follows that $P_nP_m=P_mP_n=P_n$. Clearly $\|P_nx-x\|\Lim{n\til\infty}0$ for all $x\in {\rm span}\big(\bigcup_{n\geqslant 1}F_{T_n}\big)$ and so, as the $P_n$ are uniformly bounded, this holds for all $x\in \overline{\rm span}\big(\bigcup_{n\geqslant 1}F_{T_n}\big)$.
By \cite{LT}, p. 47, it follows that the $ \dim F_{T_1}$-dimensional subspaces $(P_{n}-P_{n-1})[X]$ form a finite-dimensional decomposition of the complemented closed subspace $\overline{\rm span}\big(\bigcup_{n\geqslant 1}F_{T_n}\big)$. Since the dimensions of summands are uniformly bounded, we can further refine the decomposition to a Schauder basis for $\overline{\rm span}\big(\bigcup_{n\geqslant 1}F_{T_n}\big)$.
\end{proof}

\begin{cor} \label{finitedimensionalisometries}
Suppose $X$ is a separable,  reflexive Banach space and $G\leqslant GL_{f}(X)$ is  a bounded subgroup,  strongly closed in $GL(X)$. Then either $G$ acts nearly trivially on $X$ or $X$ has a complemented subspace with a Schauder basis.
\end{cor}

\begin{proof}
By Theorem \ref{decomp}, $X=K_G\oplus (K_{G^*})_\perp$ and $ (K_{G^*})_\perp$ is either infinite-dimensional or $X=K_G$. 
If $X=K_G$, then $G$ is almost periodic and hence acts nearly trivially on $X$ by Theorem \ref{christian}. 

On the other hand, if $(K_{G^*})_\perp$ is infinite-dimensional, 
$$
G|_{ (K_{G^*})_\perp}=\{T|_{ (K_{G^*})_\perp}\del T\in G\}
$$
is a bounded subgroup of $GL_f( (K_{G^*})_\perp)$ having no non-zero relatively compact orbits. By Theorem \ref{schauder}, $ (K_{G^*})_\perp$ and thus also $X$ has a complemented subspace with a Schauder basis.
\end{proof}

By the same reasoning, we obtain the following result.
\begin{thm}
Suppose $X$ is separable, reflexive and $G\leqslant GL_f(X)$ is a bounded subgroup. Then $G$ acts nearly trivially on $X$ or $X$ has a complemented subspace with a finite-dimensional decomposition.
\end{thm}

\begin{proof}
If $G$ does not act nearly trivially on $X$, we can choose a sequence $T_1,T_2,\ldots \in G$ such that
$$
F_{T_1}+\ldots+F_{T_n}\subsetneq F_{T_1}+\ldots+F_{T_{n+1}}
$$
for all $n\geqslant 1$. We can now repeat the proof of Theorem \ref{schauder} to get a finite-dimensional decomposition of the complemented subspace $\overline{\rm span}\big(\bigcup_{n\geqslant 1}F_{T_n}\big)$.
\end{proof}


\section{Groups of finite-dimensional isometries}
\subsection{Norm closed subgroups of $GL_f(X)$}

\begin{lemme}\label{HTU}
Let $X$ be a Banach space and $T, U\in GL_f(X)$. Then there exists an injective map from
$H_{TU}/(H_T \cap H_U)$ into $F_T \cap F_U$.
In particular, if $F_T \cap F_U=\{0\}$, then $H_{TU}=H_T \cap H_U$.
\end{lemme}
\begin{proof}
Obviously $H_T \cap H_U \subseteq H_{TU}$, and
 $U-\Id$ induces a map from $H_{TU}/(H_T \cap H_U)$ into $F_U$. For any $x \in H_{TU}$,
$$(U-\Id)x=(U-TU)x=(\Id-T)(Ux) \in F_T,$$
so the map induced by $U-\Id$ takes its values in $F_T \cap F_U$.

Finally, if $(U-\Id)x=0$ for $x \in H_{TU}$, then $Ux=x$ and $Tx=TUx=x$, so $x \in H_T \cap H_U$. This proves that the map induced by $U-\Id$ is injective.
\end{proof}

\begin{lemme}\label{triangle}
Let $X$ be a Banach space and  $T,U\in GL_f(X)$. Then
$$|\dim F_T-\dim F_U| \leqslant \dim F_{TU} \leqslant \dim F_T+\dim F_U.$$
\end{lemme}

\begin{proof} Since $TU-\Id=(T-\Id)+T(U-\Id)$, we see that $F_{TU} \subseteq F_T + T[F_U]$ and therefore
$$
\dim F_{TU} \leqslant  \dim F_T+\dim T[F_U]= \dim F_T +\dim F_U.
$$
Applying this to $TU$ and $U^{-1}$, it follows that
$$
\dim F_T \leqslant \dim F_{TU}+\dim F_{U^{-1}}=\dim F_{TU}+\dim F_U,
$$
and so $\dim F_{TU} \geqslant \dim F_T-\dim F_U$. Since also
$$
\dim F_{TU}=\dim F_{U^{-1}T^{-1}} \geqslant \dim F_{U^{-1}}-\dim F_{T^{-1}}=\dim F_U - \dim F_T,
$$
we find that $\dim F_{TU} \geqslant |\dim F_T-\dim F_U|$.
\end{proof}

\begin{lemme}\label{localmaximal}
Let $X$ be a Banach space  and $G\leqslant GL_f(X)$ be a bounded subgroup, norm-closed in $GL(X)$.
Then for any non-empty, norm-open $\ku U\subseteq G$, there 
is a smaller norm-open set $\tom \neq \ku V\subseteq \ku U$ such that $\dim F_{T}$ is constant for $T\in \ku V$.
\end{lemme}

\begin{proof}We work in the norm topology.
Note that for any $n$ the set
$$
D_n=\{A\in \ku L(X)\del \dim (\im \;A)\leqslant n\}
$$
is closed, whence 
$$
E_n=\{T\in G\del \dim F_{T}\leqslant n\}=G\cap \big(\Id+D_n\big)
$$
is closed in $G$. Now, if $\ku U\subseteq G$ is non-empty open, then as $G=\bigcup_{n\geqslant 1}E_n$, it follows from the Baire category Theorem in $\ku U$ that some $E_n \cap \ku U$ has non-empty interior. Moreover, if $n$ is minimal with this property, we see that $\ku V={\rm int}(E_n\cap \ku U)\setminus E_{n-1}\neq \tom$ and $\dim F_{T}=n$ for all $T\in \ku V$.
\end{proof}

From this we deduce the following lemma.
\begin{lemme}\label{uniform continuity}
Let $X$ be a Banach space and $G\leqslant GL_f(X)$ be a bounded subgroup, norm-closed in $GL(X)$.  Then there are  $\delta>0$ and a constant $N$ such that for all $T,U \in G$,
$$
\norm{T-U}<\delta \Rightarrow |\dim F_{T} - \dim F_{U}| \leqslant N.
$$
It follows that if $G$ is norm-compact, then 
$$
\sup \big\{\dim F_T  \del T \in G\big\}<\infty.
$$
\end{lemme}

\begin{proof} Without loss of generality, $G\leqslant {\rm Isom}_f(X)$. By Lemma \ref{localmaximal} there are $S \in G$ and $\delta>0$ such that $\dim F_T =\dim F_{S}$ whenever $\|T-S\| \leqslant \delta$. It follows that if $T,U \in G$ satisfy $\|T-U\|<\delta$ and thus $\|STU\inv-S\|<\delta$, then $\dim F_{STU^{-1}}=\dim F_{S}$ and so, by Lemma \ref{triangle},
\[\begin{split}
|\dim F_T-\dim F_U| &=|\dim F_T-\dim F_{U^{-1}}| \\
&\leqslant \dim F_{TU^{-1}}\\
&\leqslant	\dim F_{S\inv} +\dim F_{STU^{-1}} \\
&= 2\,\dim F_{S}.
\end{split}\]
Setting $N=2\,\dim F_S$, the result follows.

If now $G$ is norm-compact, it can be covered by a finite number of open balls of radius $\delta$, whereby  $\dim F_T$ is bounded above for $T \in G$.
\end{proof}

\begin{lemme}\label{preorthogonalofHT}
Let $X$ be a Banach space with a G\^ateaux differentiable norm and let $T$ be a finite-dimensional isometry on $X$. Then
$$
F_T=J[H_T]_{\perp}.
$$
\end{lemme}

\begin{proof} From Lemma \ref{facto} we easily deduce that $J[H_T] \subseteq H_{T^*}$ and therefore, as $F_T$ is closed,
$$
F_T=(F_T^\perp)_\perp=(H_{T^*})_{\perp} \subseteq J[H_T]_{\perp}.
$$
It follows by Lemma \ref{preorthogonal} (a) that
$$
X=H_T \oplus F_T \subseteq H_T \oplus J[H_T]_{\perp} \subseteq X,
$$
whence $F_T=J[H_T]_{\perp}$.
\end{proof}

\begin{thm}\label{discrete}
Let $X$ be a Banach space with a G\^ateaux differentiable norm, $G\leqslant {\rm Isom}_f(X)$ be norm-closed in $GL(X)$
 and  assume that no non-zero point of $X$ has a totally bounded $G$-orbit. Then $G$ is discrete in the norm topology and hence is locally finite, i.e., any finitely generated subgroup is finite.
\end{thm}

\begin{proof}
We work in the norm topology on $G$.
Assume towards a contradiction that $G$ is not discrete. Then for every neighbourhood $\ku W$ of $\Id$ there is a finite-dimensional isometry $U\in \ku W$ with $\dim F_U\geqslant 1$. By Lemma \ref{uniform continuity}, let $n \geqslant 1$ be  minimal such that for some neighbourhood $\ku V$ of ${\rm Id}$ in $G$ we have ${\rm dim}(F_U)\leqslant n$ for all $U\in \ku V$. Let $\ku U$ be a conjugacy invariant neighbourhood of ${\rm Id}$ in $G$ such that $\ku U^2\subseteq \ku V$ and choose an isometry $T\in \ku U$ such that ${\rm dim}(F_T)=n$.

Suppose first that $S$ is an isometry such that $S[F_T]\cap F_T=\{0\}$.
Then, since $F_{STS\inv}=S[F_T]$, Lemma \ref{HTU} implies that $H_{TSTS\inv}=H_T\cap H_{STS\inv}$. However, $T\in \ku U$ and $\ku U$ is conjugacy invariant, so $TSTS\inv\in \ku U^2
\subseteq \ku V$, and thus
\[\begin{split}
{\rm codim}(H_T)&={\rm dim }(F_T)\\
&=n\\
&\geqslant {\rm dim}(F_{TSTS\inv})\\
&={\rm codim}(H_{TSTS\inv})\\
&={\rm codim}(H_T\cap H_{STS\inv}),
\end{split}\]
whence, as ${\rm codim}(H_T)={\rm codim}(H_{STS\inv})$, we have $H_T=H_{STS\inv}=S[H_T]$.

Now, by Lemma \ref{loin},  any isometry $U$ can be written as a product of two isometries $S_1, S_2$ in $G$ such that $S_i[F_T]\cap F_T=\{0\}$. By the above calculation, $H_T$ is $S_i$-invariant and therefore also $U$-invariant. On the other hand, by Lemma \ref{facto}, we have that $U^*J[H_T]=JU^{-1}[H_T]=J[H_T]$, so, as $U$ is arbitrary, $J[H_T]$ is $G$-invariant.

Since, by Lemma \ref{preorthogonalofHT}, $F_T=J[H_T]_{\perp}$,  it  is  a non-trivial, finite-dimensional, $G$-invariant subspace of $X$, contradicting that no orbit is totally bounded.

 To see that $G$ is locally finite, note that by Lemma \ref{finitely generated} any finitely generated subgroup of ${\rm Isom}_f(X)$ is precompact in the norm topology and therefore, being discrete, must be finite.
\end{proof}

\begin{thm}\label{discretebis}
Let $X$ be a Banach space with separable dual, $G\leqslant GL_f(X)$ be a bounded subgroup,  norm-closed in  $GL(X)$ and assume that no non-zero point of $X$ has a totally bounded $G$-orbit. Then $G$ is discrete and locally finite in the norm topology.
\end{thm}

\begin{proof}
By renorming $X$, we may assume $G \leqslant {\rm Isom}_f(X)$. Thus, by Proposition \ref{lancien} and  the fact that a norm whose dual norm is LUR is G\^ateaux differentiable, we can also suppose that the norm on $X$ is G\^ateaux differentiable, so the result follows from Theorem \ref{discrete}.
\end{proof}


\subsection{Decompositions by strongly closed subgroups of $GL_f(X)$.}
\begin{thm}\label{decompos}
Suppose $X$ is a separable, reflexive Banach space and $G\leqslant GL_f(X)$ is bounded and strongly closed in $GL(X)$. 
Let also $G_0\leqslant G$ denote the connected component of the identity in $G$.

Then $G_0$ acts nearly trivially on $X$ and therefore is a compact Lie group. Moreover, $G_0$ is open in $G$, while $G/G_0$ is a countable, locally finite group. It follows that $G$ is an amenable Lie group.

Furthermore, $X$ admits a $G$-invariant decomposition $X=X_1\oplus X_2\oplus X_3\oplus X_4$, where
\begin{enumerate}
\item no non-zero point of $X_1$ has a relatively compact $G$-orbit,
\item every $G$-orbit on $X_2\oplus X_3$ is relatively compact,
\item $X_4$ is the subspace of points which are fixed by $G$,
\item $X_2$ is finite-dimensional and $X_1\oplus X_3\oplus X_4$ is the subspace of points which are fixed by $G_0$,
\item if $X_1\neq\{0\}$, then $X_1$ has a complemented subspace with a Schauder basis, while if $X_1=\{0\}$, then  $G$ acts nearly trivially on $X$.
\end{enumerate}
\end{thm}

\begin{proof}
Note that, by Proposition \ref{topos}, since $X^*$ is separable and $G\leqslant GL_f(X)$, the norm and strong operator topologies coincide on $G$. 

Let now $X=K_G\oplus (K_{G^*})_\perp$ be the $G$-invariant decomposition given by Theorem \ref{decomp} and, for simplicity of notation, let $X_1=(K_{G^*})_\perp$. 

We claim that $G|_{X_1}=\{T|_{X_1}\del T\in G\}$ is a strongly closed subgroup of $GL(X_1)$.  To see this,
suppose $T_n\in G$ and $T\in GL(X_1)$ are such that $T_n|_{X_1}\Lim{\text{SOT}}T$. 
By Proposition \ref{Gcompact}, $G|_{K_G}=\{T|_{K_G}\del T\in G\}\leqslant GL(K_G)$ is precompact in the strong operator topology and thus $({T_n}|_{K_G})_{n\in \N}$ has a subsequence $(T_{n_i}|_{K_G})_{i\in \N}$ converging in the strong operator topology to an operator $S\in GL(K_G)$, whence $T_{n_i}\Lim{\text{SOT}}S\oplus T\in GL(K_G\oplus X_1)$. 
Since $G$ is closed in the strong operator topology, we see that $S\oplus T\in G$, whence $T=(S\oplus T)|_{X_1}\in G|_{X_1}$, showing that $G|_{X_1}$ is strongly closed. 

Now, by Proposition \ref{topos}, since $X_1^*$ is separable and $G|_{X_1} \leqslant GL_f(X_1)$, the norm and strong operator topologies coincide on $G|_{X_1}$. Moreover, as there is no non-zero relatively compact $G|_{X_1}$-orbit on ${X_1}$, Theorem \ref{discretebis} implies that $G|_{X_1}$ is discrete and locally finite. Being separable, it must be countable. Moreover, since the map $ T\in G\mapsto T|_{X_1}\in G|_{X_1}$ is strongly continuous, its kernel $G'=\{T\in G\del {X_1}\subseteq H_T\}$ is a countable index, closed normal subgroup of the Polish group $G$ and thus must be open by the Baire category Theorem. Since every $G'$-orbit on $X$ is relatively compact, Proposition \ref{Gcompact} implies that $G'=\overline{G'}^{\text{SOT}}$ is compact in the strong operator topology and so, by Theorem \ref{christian}, 
$G'$ acts nearly trivially on $X$, and since $G_0 \leqslant G'$, the same holds for $G_0$. It follows that $G'$ is a compact Lie group, whence also $G$ is Lie, and that $G_0$ is a compact, open, normal subgroup in $G$. Since $G'$ is compact, $[G'\colon G_0]<\infty$, and so $G/G_0$ is an extension of the finite group $G'/G_0$ by the locally finite group $G/G'$ and thus is itself locally finite.
Thus, as both $G_0$ and $G/G_0$ are amenable, so is $G$.

We now consider the canonical complement $K_G=H_G\oplus \big((H_{G^*})_\perp\cap K_G\big)$ to $X_1=(K_{G^*})_\perp$.
Set $Y=\big((H_{G^*})_\perp\cap K_G\big)$ and let $\pi\colon G\til GL(Y)$ be the strongly continuous representation $\pi(T)=T|_Y$. Since $T|_{X_1}=\Id_{X_1}$ for all $T \in G_0$, $\pi(G_0)$ is strongly closed in $GL(Y)$, and so, by Theorem \ref{christian}, $Y$ admits a $G_0$-invariant decomposition 
$$
Y=H_{\pi(G_0)}\oplus (H_{\pi(G_0)^*})_\perp,
$$
with $X_2=(H_{\pi(G_0)^*})_\perp$ finite-dimensional. Moreover, since  $G_0$ is normal in $G$, both $X_3=H_{\pi(G_0)}$ and $(H_{\pi(G_0)^*})_\perp$ are $G$-invariant. So, letting $X_4=H_G$, we obtain the desired decomposition. 

If $X_1 \neq \{0\}$, then the existence of the complemented subspace of $X_1$ with a Schauder basis follows from Theorem \ref{schauder}. 
If $X_1=\{0\}$, then  $G$ is almost periodic and therefore acts nearly trivially on $X$ by Theorem  \ref{christian}. 
\end{proof}


\section{Spaces with a small algebra of operators}

\subsection{Unconditional sequences of eigenvectors}
\begin{lemme} Let $0<\theta<\epsilon<\pi$ and $0\leqslant \alpha \leqslant 2\pi$. For any interval  $I\subseteq \N$ of cardinality at least $3\pi/\theta+1$, there exists a subinterval $J\subseteq I$ of cardinality at least $\epsilon/\theta$ such that for all $m \in J$,  $d(m\theta,\alpha+2\pi\N) <\epsilon$.
\end{lemme}

\begin{proof} Let $(u_n)_{n \in I}$ be the finite sequence $(n\theta-\alpha)_{n \in I}$ and note that the difference between the first and last element of $u_n$ is at least $3\pi$, while the difference between two successive elements is at most $\theta<\pi$. Write the first element of the sequence as $k2\pi+\alpha_0$, for some $\alpha_0 \in [0,2\pi[$, and let $\alpha'=k2\pi+\alpha$ if $\alpha_0 \leqslant \alpha$, or $\alpha'=(k+1)2\pi+\alpha$ otherwise.

Let $n_1$ be the largest element of $I$ such that $u_{n_1} \leqslant \alpha'$, and note that the difference between $u_{n_1}$ and the first element of the sequence is at most $2\pi$. Now consider  $p \in \N$ maximal  with $1 \leqslant p < \epsilon/\theta$ and note that $u_{n_1+p}-u_{n_1} \leqslant p\theta <\epsilon$. Therefore since $\alpha' \in [u_{n_1},u_{n_1+p}]$, $d(u_n,\alpha')=d(n\theta,\alpha+2\pi\N)<\epsilon$ for any $n$ in the interval $J=\{n_1,\ldots,n_1+p\}$. This interval has cardinality at least $\epsilon/\theta$. The difference between the first $u_n$ and $u_{n_1+p}$ is at most $2\pi+p\theta \leqslant 2\pi + \epsilon<3\pi$, therefore  $J$ is a subinterval of $I$.
\end{proof}

\begin{lemme}\label{roundtheclock} 
Let $(\epsilon_k)_{k \geqslant 1}$ be a sequence of real numbers in the interval $]0,\pi[$. Let $(\theta_k)_{k \geqslant 1}$ be a  sequence of positive real numbers such that
$\epsilon_1/\theta_1 \geqslant 1$ and $\epsilon_k/\theta_k>3\pi/\theta_{k-1}+1$ for all $k>1$.
Let  $n \in \N$ and $(\alpha_k)_{1 \leqslant k \leqslant n}$ be a finite sequence of real numbers.
Then there exists some  $m \in \N$ such that for all $k=1,\ldots, n$
$$
d(m\theta_k,\alpha_k+2\pi\N)<\epsilon_k.
$$
\end{lemme}

\begin{proof} First we may by the previous lemma pick an interval $J_n$ of cardinality at least
$\epsilon_n/\theta_n$ such that for all $m \in J_n$,
$d(m\theta_n-\alpha_n,2\pi\N)<\epsilon_n$.
Since $\epsilon_n/\theta_n-1>3\pi/\theta_{n-1}$, the lemma applies again to obtain an interval $J_{n-1} \subseteq J_n$ of cardinality at least
$\epsilon_{n-1}/\theta_{n-1}$ such that for all $m \in J_{n-1}$,
$d(m\theta_{n-1}-\alpha_{n-1},2\pi\N)<\epsilon_{n-1}$.
After $n$ steps we have obtained the desired result for any $m$ in some interval $J_1$ of cardinality at least $\epsilon_1/\theta_1$ and therefore non-empty by hypothesis.
\end{proof}

\begin{prop}\label{eigen}Let $T$ be a linear contraction on a complex Banach space $X$ having infinitely many eigenvalues of modulus $1$. Then, for any $\eps>0$, $X$ contains a $(1+\eps)$-unconditional basic sequence.  
\end{prop}

\begin{proof}Modulo replacing $T$ by some $\lambda T$, $|\lambda|=1$, we can suppose there are distinct $\theta_n$ converging to $0$ such that each $e^{i\theta_n}$ is an eigenvalue for $T$ with corresponding normalised eigenvector $x_n\in X$. Moreover, without loss of generality, we can also assume  that $0<\theta_n <\pi$ and that, given $\epsilon_n>0$ such that $\sum_{n=1}^\infty \epsilon_n\leqslant \delta$, where $\frac{2+\delta}{2-\delta}<1+\eps$, one has 
$\epsilon_1/\theta_1 \geqslant 1$ and $\epsilon_n/\theta_n>3\pi/\theta_{n-1}+1$ for all $n>1$.

Suppose $N \geqslant 1$ and consider a vector of the form $x=\sum_{n \leqslant N}\lambda_n x_n$ for $\lambda_n\in \C$.
Now, if $e^{i\alpha_n}\in \T$ and $m\geqslant 1$, we have 
$$
T^m(x)=\sum_{n \leqslant N}e^{im\theta_n}\lambda_n x_n,
$$  
and so
$$ 
\Norm{T^m(x) -\sum_{n \leqslant N}e^{i\alpha_n}\lambda_n  x_n}=\Norm{{\sum_{n \leqslant N}( e^{im\theta_n}-e^{i\alpha_n})\lambda_n x_n}} \leqslant \sum_{n \leqslant N}|1-e^{i(\alpha_n-m\theta_n)}||\lambda_n|.
$$
Choosing $m$ according to Lemma \ref{roundtheclock} so that $|\alpha_n-m\theta_n|$ is sufficiently small for each $n\leqslant N$,  we conclude that 
\begin{equation}\label{equa1}
\Norm{T^m(x)-\sum_{n \leqslant N} e^{i\alpha_n}  \lambda_n x_n} 
\leqslant
\sum_{n \leqslant N}\epsilon_n|\lambda_n|
\leqslant \delta\cdot \sup_{n\leqslant N}|\lambda_n|
\end{equation}
and so also
\begin{equation}\label{equa2}
\Norm{\sum_{n \leqslant N} e^{i\alpha_n} \lambda_n x_n} \leqslant \Norm{\sum_{n \leqslant N} \lambda_n x_n} +
 \delta\cdot \sup_{n\leqslant N}|\lambda_n|. 
\end{equation}
Let $n_0 \leqslant N$ be chosen such that $|\lambda_{n_0}|=\sup_{n\leqslant N}|\lambda_n|$ and set 
$e^{i\alpha_{n_0}}=1$ and $e^{i\alpha_n}=-1$ for $n\neq n_0$. By (2), we then have 
\[\begin{split}
2|\lambda_{n_0}| 
&\leqslant 
\Norm{\lambda_{n_0}x_{n_0}+\sum_{\substack{^{n \leqslant N}_{ n \neq n_0}}}\lambda_n x_n}
+\Norm{\lambda_{n_0}x_{n_0}-\sum_{\substack{^{n \leqslant N}_{ n \neq n_0}}}\lambda_n x_n}\\
&\leqslant 
2\norm{\sum_{n \leqslant N}\lambda_n x_n} + \delta\cdot |\lambda_{n_0}|,
\end{split}\]
i.e., $\sup_{n\leqslant N}|\lambda_n|=|\lambda_{n_0}|\leqslant \frac{2}{2-\delta}\norm{\sum_{n \leqslant N}\lambda_n x_n}$. Combining this with (2), for any $\lambda_n$ and $e^{i\alpha_n}\in \T$, we have
$$
\Norm{\sum_{n \leqslant N} e^{i\alpha_n} \lambda_n x_n} \leqslant \frac{2+\delta}{2-\delta}\Norm{\sum_{n \leqslant N} \lambda_n x_n},
$$
which shows that $(x_n)_{n=1}^\infty$ is a $(1+\eps)$-unconditional basic sequence.
\end{proof}

With this in hand we may prove the following:

\begin{thm}\label{posfin} 
Let $X$ be a Banach space containing no unconditional basic sequence. Then any bounded subgroup $G\leqslant GL_{in}(X)$ is contained in $GL_f(X)$ and so, in particular, 
$$
{\rm Isom}_{af}(X)={\rm Isom}_f(X).
$$
\end{thm}

\begin{proof}
By Theorem \ref{in then af}, it suffices to prove that any bounded subgroup $G\leqslant GL_{af}(X)$ is contained in $GL_f(X)$ and, by renorming, that any almost finite-dimensional isometry is finite-dimensional.

Assume first that $X$ is complex.  Since almost finite-rank operators are Riesz operators, the spectrum of an almost finite-dimensional isometry $T$ of $X$ is either a finite or an infinite sequence of distinct eigenvalues with finite multiplicity, together with the value $1$, which is the limit of the sequence in the infinite case  \cite{GM}. Since $X$ does not contain an unconditional basic sequence, by Proposition \ref{eigen}, the second case cannot occur and so $\sigma(T)$ is finite. Lemma \ref{gelfand} then implies that $T-\Id$ must have finite rank.

If instead $X$ is real, we claim that the complexification $\hat{X}$ of $X$ does not contain an unconditional basic sequence either. For if it did, then $X \oplus X$ would in particular contain a real unconditional sequence, spanning a real subspace $Y$.  And since some subspace of $Y$ must either embed into the first or second summand of the decomposition $X \oplus X$, $X$ itself would contain a real unconditional basic sequence, contradicting our assumption.

Now, as mentioned in Section \ref{complexification}, if $T$ is an almost finite-dimensional isometry of $X$, $\hat{T}$ is an almost finite-dimensional isometry of $\hat{X}$, which then must be finite-dimensional. Since $F_{\hat{T}}=F_T+iF_T$, we finally conclude that also $F_T$ is finite-dimensional.
\end{proof}

Combining Theorem \ref{discretebis} and Theorem \ref{posfin}, we also obtain

\begin{cor}\label{totally unbounded}
Let $X$ be a Banach space with separable dual and not containing an unconditional basic sequence. Let $G \leqslant GL_{af}(X)$ be a bounded subgroup, and assume that no non-zero point of $X$ has a totally bounded $G$-orbit. Then $G$ is discrete and locally finite in the norm topology.
\end{cor}

\begin{proof} Since $GL_{af}(X)$ is norm-closed in $GL(X)$, by replacing $G$ by $\overline{G}^{\|\cdot\|}$ we may assume that $G$ is norm-closed. Also, by renorming $X$, we can suppose that $G\leqslant {\rm Isom}_{af}(X)$. Then, since by Theorem \ref{posfin}, $G \leqslant {\rm Isom}_f(X)$, Theorem \ref{discretebis} applies.
\end{proof}


\subsection{Groups of isometries in spaces with few operators}
We shall now combine the results of the previous sections in order to obtain a description of the group of isometries of  a separable reflexive space $X$ with a small algebra of operators. In this case, Proposition \ref{topos} will then ensure that the norm and the strong operator topology coincide on ${\rm Isom}(X)$. 

We recall that an infinite-dimensional Banach space is {\em decomposable} if it can be written as the direct sum of two infinite-dimensional subspaces and  {\em hereditarily indecomposable} (or HI) if it has no  decomposable subspace. The first construction of an HI space was given by W. T. Gowers and B. Maurey in \cite{GM} as an example of an infinite-dimensional Banach space not containing an unconditional basic sequence, since it is clear that an HI space cannot contain any subspace with an unconditional basis. Furthermore, Gowers and Maurey proved that any operator on a complex HI space is a  strictly singular perturbation of a multiple of the identity. We shall call this latter property of a (real or complex) Banach space the {\em $\lambda \Id+S$-property}. Note that the  $\lambda \Id+S$-property easily implies that the space is indecomposable, since no projection with infinite-dimensional range and corange is of the form  $\lambda \Id+S$, with $S$ strictly singular.

As an immediate consequence of Theorem \ref{posfin}, we have the following.

\begin{thm}\label{nearly trivial}
Let $X$ be a Banach space with the $\lambda\Id+S$-property and containing no unconditional basic sequence. Then each individual isometry acts nearly trivially on $X$.
\end{thm}

\begin{cor}\label{real} 
Let $X$ be a separable reflexive space with the  $\lambda \Id+S$-property and containing no unconditional basic sequence. Assume $X$ does not have a Schauder basis. Then ${\rm Isom}(X)$ acts nearly trivially on $X$.
\end{cor}

\begin{proof}
Note that, depending on whether $X$ is complex or real, we have by Theorem \ref{nearly trivial} that
$$
{\rm Isom}(X)=\T \times {\rm Isom}_f(X)
$$
or
$$
{\rm Isom}(X)=\{-1,1\} \times {\rm Isom}_f(X).
$$
Thus, by Proposition \ref{topos}, the norm and strong operator topologies coincide on ${\rm Isom}(X)$. So, as ${\rm Isom}_f(X)$ is norm closed in ${\rm Isom}(X)$, it is strongly closed in ${\rm Isom}(X)$ and hence also in $GL(X)$. By Corollary \ref{finitedimensionalisometries}, we see that  ${\rm Isom}_f(X)$ and thus also  ${\rm Isom}(X)$ acts nearly trivially on $X$.
\end{proof}

\begin{cor}\label{noschauder} Let  $X$ be a separable, reflexive, complex HI space without a Schauder basis. Then ${\rm Isom}(X)$ acts nearly trivially on $X$. 
\end{cor}

Of course, in the case when $X$ is permitted to have a Schauder basis, Theorem \ref{nearly trivial} does not itself describe the global action of the group ${\rm Isom}(X)$, but with some extra hypotheses, we can get more information. Observe first that if $X$ is an infinite-dimensional space with the $\lambda\Id+S$-property, then $GL(X)=\K^\star\times GL_s(X)$.

\begin{thm}
Let $X$ be a separable, infinite-dimensional, reflexive space with the $\lambda \Id +S$-property and containing no unconditional basic sequence. Assume $G \leqslant GL_s(X)$ is a bounded subgroup. Then either
\begin{itemize}
\item[(i)] $G$ acts nearly trivially on $X$, or
\item[(ii)] $X$ admits a $G$-invariant decomposition $F \oplus  H$, where $F$ is finite-dimensional, and $G|_H=\{T|_H\del T\in G\}$ is a countable, discrete, locally finite subgroup of $GL_f(H)$ none of whose non-zero orbits are relatively compact. Moreover, in this case, $X$ has a Schauder basis.
\end{itemize}
\end{thm}

\begin{proof}By Theorem \ref{posfin}, $G\leqslant GL_f(X)$. Also, as in the proof of Corollary \ref{real}, we see that  ${\rm Isom}_f(X)$ is strongly closed in $GL(X)$.
By renorming, we may assume that $G \leqslant {\rm Isom}_{f}(X)$, whence also $M=\overline{G}^{\rm SOT}$ is a subgroup of ${\rm Isom}_f(X)$. Applying Theorem \ref{decompos} to $M$, we obtain a decomposition 
$$
X=X_1\oplus X_2\oplus X_3\oplus X_4.
$$ 

Here either $X_1=\{0\}$, in which case (i) holds, or $X_1$ has a complemented subspace with a Schauder basis. 

By the indecomposability of $X$, in the latter case, $X$ itself has a Schauder basis, $H=X_1$ has finite codimension in $X$ and
$F=X_2 \oplus X_3 \oplus X_4$ is finite-dimensional. 
Furthermore, $M|_{X_1}=\{T|_{X_1}\del T\in M\}$ is a strongly closed subgroup of $GL(X_1)$ contained in ${\rm Isom}_f(X_1)$ such that no non-zero $M$-orbit on $X_1$ is relatively compact. 
 By Proposition \ref{topos}, the strong operator and norm topologies coincide on $M|_{X_1}$, while, by Theorem \ref{discretebis}, $M|_{X_1}$ is countable, discrete and locally finite. It follows that $G|_{X_1}=M|_{X_1}$.
 \end{proof}

While the real HI space constructed by Gowers and Maurey has the $\lambda \Id +S$-property, this does not generalize to all real HI spaces \cite{F2}. Nor is it true that the complexification of a real HI space is complex HI, see \cite{F2} and Proposition 3.16 \cite{GH}. So it is not clear whether Corollary \ref{noschauder} extends to the real case.

Our methods may also be used in spaces with the $\lambda \Id+S$-property containing unconditional basic sequences. Such spaces do exist, see, for example, \cite{AF}. 

\begin{thm} 
Let $X$ be a separable reflexive space with the $\lambda \Id+S$-property. Then either ${\rm Isom}(X)$ acts nearly trivially on $X$, or $X$ admits an isometry invariant decomposition $X=F \oplus Y$, where $F$ is finite-dimensional and no orbit of a non-zero point of $Y$ under ${\rm Isom}(X)$ is totally bounded.
\end{thm}

\begin{proof} 
By Theorem \ref{in then af}, depending on whether $X$ is complex or real, we have that
$$
{\rm Isom}(X)=\T \times {\rm Isom}_{af}(X)
$$
or
$$
{\rm Isom}(X)=\{-1,1\} \times {\rm Isom}_{af}(X),
$$
and so Proposition \ref{topos} applies to deduce that the strong operator and the norm topologies coincide on ${\rm Isom}(X)$, whence $G={\rm Isom}_{af}(X)$ is strongly closed in $GL(X)$.

By Theorem \ref{decomp}, we have a $G$-invariant decomposition
$X=K_G \oplus (K_G^*)^{\perp}$, where either $X=K_G$ or $K_G$ has infinite codimension. Furthermore, $X$ is indecomposable, so either $X=K_G$ or $K_G$ is finite-dimensional.
In the first case all $G$ is almost periodic and thus  by Theorem \ref{christian}  acts nearly trivially on $X$.
In the second case, we define $F=K_G$ and $Y=(K_G^*)^{\perp}$.
\end{proof}


\section{Maximality and transitivity in spaces with few operators}

Let us begin by reviewing the various types of norms defined and studied by Pe\l czy\'nski and Rolevicz in \cite{PR,R}.

\begin{defi} Let $(X,\norm\cdot)$ be a Banach space and for any $x \in S_X$ let ${\mathcal O}(x)$ denote the orbit of $x$ under the action of ${\rm Isom}(X,\norm\cdot)$. The norm $\norm\cdot$ on $X$ is
\begin{itemize}
\item[(i)] {\em transitive} if for any $x \in S_X$, ${\mathcal O}(x)=S_X$,
\item[(ii)] {\em almost transitive} if for any $x \in S_X$, ${\mathcal O}(x)$ is dense in $S_X$,
\item[(iii)] {\em convex transitive} if for any $x \in  S_X$, ${\rm conv}\,{\mathcal O}(x)$ is dense in $B_X$,
\item[(iv)] {\em uniquely maximal} if whenever $\triple\cdot$ is an equivalent norm on $X$ such that
${\rm Isom}(X,\norm\cdot) \leqslant {\rm Isom}(X,\triple\cdot)$, then $\triple\cdot$ is a scalar multiple of $\norm\cdot$.
\item[(v)] {\em maximal} if whenever $\triple\cdot$ is an equivalent norm on $X$ such that
${\rm Isom}(X,\norm\cdot) \leqslant {\rm Isom}(X,\triple\cdot)$, then ${\rm Isom}(X,\norm\cdot)={\rm Isom}(X,\triple\cdot)$.
\end{itemize}
\end{defi}

Here, the implications (i)$\Rightarrow$(ii)$\Rightarrow$(iii) as well as (iv)$\Rightarrow$(v) are obvious. Furthermore, Rolewicz \cite{R} proved that any convex transitive norm must be uniquely maximal and later  E. R. Cowie \cite{C} reversed this implication by showing that a uniquely maximal norm is convex transitive. So (i)$\Rightarrow$(ii)$\Rightarrow$(iii)$\Leftrightarrow$(iv)$\Rightarrow$(v).

Almost transitive norms are not too difficult to obtain. For example, the classical norm on  $L_p([0,1])$, $1 \leqslant p<\infty$, is almost transitive \cite{R,GJK}. It is also known that the non-trivial ultrapower of a space with an almost transitive norm will have a transitive norm. There are therefore many examples of non-separable, non-Hilbertian spaces with a transitive norm.

Rolewicz \cite{R} also proved that if a space has a $1$-symmetric basic sequence, then the norm is maximal. Therefore the usual norms on the spaces $c_0$ and $\ell_p$, $1 \leqslant p<\infty$, are maximal, though they are not convex transitive.

More interesting than asking whether a specific norm on a Banach space $X$ has one of the above forms of transitivity or maximality is the question of whether $X$ admits an equivalent norm with these properties. In this direction,  J. Becerra Guerrero and  A. Rodr\'iguez-Palacios \cite{becerra} showed the following interesting fact.

\begin{thm}[Becerra Guerrero and  Rodr\'iguez-Palacios]\label{conv trans}
Suppose $X$ is either Asplund or has the Radon--Nykodym property and that the norm $\norm\cdot$ of $X$ is convex transitive. Then $\norm\cdot$ is almost transitive, uniformly convex and uniformly smooth. In particular, $X$ is super-reflexive.
\end{thm}

This gives a  list of spaces with no equivalent convex transitive norm, $c_0$, $\ell_1$, Tsirelson's space $T$, Schlumprecht's space $S$, and Gowers-Maurey's space $GM$, for example.

As another application of Theorem \ref{conv trans}, we see that Theorem \ref{schauder} may be improved when we assume that the norm on $X$ is convex transitive.

\begin{thm}\label{convextransitive} 
Let $X$ be a separable, reflexive, Banach space with a convex transitive norm that carries a non-trivial finite-dimensional isometry. Then $X$ has a Schauder basis.
\end{thm}

\begin{proof}
Note that by Theorem \ref{conv trans}, the norm is almost transitive. So we can repeat the proof of Theorem \ref{schauder} making sure that the $S_2,S_3,\ldots\in {\rm Isom}(X)$ are chosen such that $X=\overline{\rm span}\big(\bigcup_{n\geqslant 1}F_{T_n}\big)$. 
\end{proof}

W. Lusky \cite{Lu} proved that every separable Banach space $Y$ is $1$-complemented in some almost transitive separable space. So there are many different spaces with almost transitive norms, and, depending on the choice of $Y$, examples without a Schauder basis. Note, however, that Lusky's theorem cannot be improved to include reflexivity, that is, if $Y$ is reflexive, but not super-reflexive, then, by Theorem \ref{conv trans}, $Y$ does not even embed into  a reflexive space with a convex transitive norm.

However, the question of whether any super-reflexive space admits an equivalent almost transitive norm has remained open hitherto. This question is due to  Deville, Godefroy, and Zizler \cite{DGZ}. Based on results by C. Finet \cite{Fi},  they observed that a positive answer  would imply that if a Banach space $X$ has an equivalent norm with modulus of  convexity of type $p \geqslant 2$ and another equivalent norm with modulus of smoothness of type $1 \leqslant q \leqslant 2$, then $X$ has an equivalent norm with both of these properties, which would be exceedingly useful in renorming theory. However, as we shall see, this approach does not work.

We now answer the question of Deville, Godefroy and Zizler along with Wood's problems of whether any Banach space has an equivalent maximal norm or even whether any bounded subgroup $G\leqslant GL(X)$ is contained in a maximal bounded subgroup.

\begin{lemme}\label{valentin}
There is a separable, infinite-dimensional, super-reflexive, complex HI space without a Schauder basis.
\end{lemme}

\begin{proof}
Recall that for a space $X$, if we define
$$
p(X)=\sup\{t\leqslant 2\del  X\ {\rm\ has\ type\ }  t\},$$
and
 $$
 q(X)=\inf\{c\geqslant 2\del  X\ {\rm\ has\ cotype\ } c\},
 $$ 
then $X$ is said to be {\em near-Hilbert} if $p(X)=q(X)=2$. Now, by \cite{F}, there exists a uniformly convex and therefore super-reflexive complex HI space $X$ with a Schauder basis.  Moreover, by \cite{F}, for any $1<p<2$,  $X$ may be chosen so that
for any finite sequence $x_1,\ldots,x_n$ of successive normalized vectors on the basis of $X$, 
we have 
$$
\norm{x_1+\cdots+x_n} \geqslant \frac{n^{1/p}}{\log_2(n+1)},
$$
which implies that $X$ does not have type more than $p$. In particular, $X$ is not near-Hilbert, and
by  classical results of A. Szankowski \cite{Sz}, see \cite{LT} Vol. II, Theorem 1.g.6, this implies
that some subspace $Y$ of $X$ does not have the Approximation Property and therefore fails to have  a Schauder basis. On the other hand, $Y$ is still super-reflexive and HI.
\end{proof}

More precise estimates about type and cotype are given in \cite{CFM} and imply that for any choice of parameters in the construction of \cite{F}, the space $X$ is not near-Hilbert and therefore has a subspace without the Approximation Property. On the other hand, it is also proved in \cite{CFM} that for any $\epsilon>0$, the space $X$ (and therefore $Y$ as well) may be chosen to have  type $2-\epsilon$ and cotype $2+\epsilon$.

\begin{prop}\label{corollaire}
Let $X$ be an infinite-dimensional Banach space with norm $\norm\cdot$. Assume that ${\rm Isom}(X,\norm\cdot)$ acts nearly trivially on $X$.  Then
$\norm\cdot$ is not maximal.
\end{prop}

\begin{proof}
Let $\norm\cdot$ be such a norm on $X$ and let $X=F \oplus H$ be the associated decomposition for which $F$ is finite-dimensional and ${\rm Isom}(X)$ acts trivially on $H$.
We fix a norm $1$ vector $x_0$ in $H$, write $H$ as a direct sum $H=[x_0] \oplus M$ and
define an equivalent norm $\triple\cdot$ on $X$
by
$$
\triple{f+\alpha x_0 +m}=\|f\|+|\alpha|+\|m\|,
$$
when $f \in F$, $\alpha$ is scalar, and $m \in M$.

If $T$ is an isometry for $\norm\cdot$, then $T|_{H}=\lambda \Id_H$ with $|\lambda|=1$, and so   the equalities
$$
\triple{T(f+\alpha x_0+m)}=\triple{Tf+\alpha\lambda x_0+\lambda m}=\|Tf\|+|\alpha|+\|m\|$$
and
$$
\triple{f+\alpha x_0+m}=\|f\|+|\alpha|+\|m\|$$
show that $T$ is also a $\triple\cdot$-isometry. Therefore, any $\norm\cdot$-isometry is a $\triple\cdot$-isometry.

Furthermore, the map $L$ on $X$ defined by
$$L(f+\alpha x_0+m)=f-\alpha x_0+m$$
is a $\triple\cdot$-isometry but not a $\norm\cdot$-isometry, since there is no scalar $\lambda$ for which $L|_{H} =\lambda \Id_H$.
This shows that $\norm\cdot$ is not  maximal. 
\end{proof}

\begin{thm}\label{nomaximal} 
There exists a separable, super-reflexive, complex Banach space $X$ that  admits no equivalent
maximal  norm. In fact, if $\norm\cdot$ is any equivalent norm on $X$, then ${\rm Isom}(X,\norm\cdot)$ acts nearly trivially on $X$.
\end{thm}

\begin{proof}
Let $X$ be the space given by Lemma \ref{valentin} and notice that, by Corollary \ref{noschauder}, the isometry group acts nearly trivially on $X$ for any equivalent norm. In particular, by Proposition \ref{corollaire}, $X$ cannot have an equivalent maximal norm.
\end{proof}

Observe that in our example, every orbit under the group of isometries is compact. In this sense, it is a particularly strong counterexample to the question of Deville, Godefroy and Zizler, so one may suspect there to be weaker counter-examples among more classical spaces.

In the real case we obtain the following counter-example to Wood's questions, which, however, is not super-reflexive.

\begin{thm}\label{realbis} There exists a real, separable, reflexive Banach space with no equivalent maximal norm.
\end{thm}

\begin{proof} The reflexive HI space $GM$ of Gowers and Maurey \cite{GM}   does not have type $p>1$. Therefore the results of Szankowski \cite{Sz} imply that some subspace $Y$ of $GM$ fails to have a Schauder basis. Furthermore,  as $Y$ contains no unconditional basic sequence and since by \cite{GM} every operator from a subspace of $X$ into $X$ is a strictly singular perturbation of a multiple of the inclusion map, $Y$ satisfies the $\lambda \Id+S$-property. The result then follows from
Theorem \ref{real} and Proposition \ref{corollaire}.
\end{proof}

In connection with this, we should mention the following conjecture due to K. Jarosz \cite{J}.

\begin{conj}[K. Jarosz] Suppose $G$  is a group and $X$ is a real Banach space with $\dim X \geqslant |G|$. Then $X$ admits an equivalent renorming such that ${\rm Isom}(X)\iso\{-1,1\} \times G$.
\end{conj}

In \cite{FG} this was verified for finite groups $G$ and separable spaces $X$, but, as we shall see, the conjecture fails for infinite $G$.
\begin{prop} 
If $X$ is the real HI space considered in Theorem \ref{realbis}, then for any norm on $X$, the group of isometries on $X$ is either finite or of cardinality $2^{\aleph_0}$. In particular, Jarosz's conjecture does not hold in general.
\end{prop}

\begin{proof}
Assume that  $X$ is the real HI space considered in Theorem \ref{realbis}. Then ${\rm Isom}(X)$ acts nearly trivially on $X$ and hence is a compact Lie group. It follows that ${\rm Isom}(X)$ is either finite or of size $2^{\aleph_0}$.
\end{proof}


\section{Questions and further comments}
\subsection{Problems concerning spaces with few operators}

We suspect that the usage of reflexivity and the non-existence of a Schauder basis in Theorem \ref{nomaximal}  are not necessary. That is, we conjecture the following.

\begin{conj} \label{conj HI}
Let $X$ be a complex HI space. Then the group of isometries acts nearly trivially on $X$.
\end{conj}

A few comments on this conjecture are in order. First of all, the following result was essentially proved by Cabello-S\'anchez \cite{CS2}. As we shall state a slightly more general result than in \cite{CS2}, we give the proof of the theorem for the sake of completeness.

\begin{thm}[Cabello-S\'anchez] 
Let $G=\T \times {\rm Isom}_f(X)$, respectively $G=\{-1,1\} \times {\rm Isom}_f(X)$, be the group of nearly trivial isometries of a space $X$. If the action of $G$ on $X$ is transitive, then the norm of $X$ is euclidean, i.e., $X$ is isometric to Hilbert space. 
\end{thm}

\begin{proof} 
By transitivity, the only ${\rm Isom}_f(X)$-invariant subspaces of $X$ are the trivial ones. Therefore, by \cite{CS2} Lemma 2, there exists an ${\rm Isom}_f(X)$-invariant inner product $\langle\cdot,\cdot\rangle$ on $X$ and some $x_0 \in S_X$ such that $\langle x_0,x_0\rangle=1$. Now, if $T\in G$ and $\lambda \in \T$ is chosen such that $U=\lambda^{-1}T\in {\rm Isom}_f(X)$, then
$$
\langle Tx_0,Tx_0\rangle =\langle\lambda Ux_0,\lambda Ux_0\rangle =|\lambda|^2\langle Ux_0,Ux_0\rangle=1
$$
as well. By transitivity of $G$, this implies that $\langle x,x\rangle=\norm{x}^2$ for all $x \in X$, which proves the theorem. \end{proof}

So by Theorem \ref{nearly trivial}, we have the following.

\begin{cor}\label{notransitive}
Let $X$ be a  space with the $\lambda \Id+S$-property and without an unconditional basic sequence. Then no equivalent norm on $X$ is transitive.
\end{cor}

As a tool towards proving Conjecture \ref{conj HI}, it may be interesting to observe that the spectrum $\sigma(T)$ depends continuously on $T$ when $T$ belongs to ${\rm Isom}_f(X)$. For, this when $X$ a complex space and $T \in{\rm Isom}_f(X)$, we denote by $F_{\lambda}(T)$ the image of the spectral projection associated to the eigenvalue $\lambda$, and note that $F_{\lambda}(T)=\ker(T-\lambda\Id)$ by Gelfand's theorem. So, by Lemma \ref{single generator}, $F_T=\bigoplus_{\lambda \neq 1} F_{\lambda}(T)$ and $H_T=F_1(T)$.

\begin{lemme}\label{norm projection}
Suppose $X$ is a complex space and $T \in {\rm Isom}_f(X)$. Then if $P_i$ denotes the canonical projection of $X$ onto $F_{\lambda_i}(T)$ corresponding to the decomposition
$$
X=F_{\lambda_1}(T)\oplus\ldots\oplus F_{\lambda_m}(T),
$$
we have $\|P_i\|=1$ for every $i$.
Since $F_T$ is the complement of $H_T=F_1(T)$, it follows that the projection onto $F_T$ has norm $\leqslant 2$.
\end{lemme}

\begin{proof}
Let $p_n(t)=n^{-1}(1+t+\cdots+t^{n-1})$. Since the operator $p_n(\lambda_i\inv T)$ has norm at most $1$ and acts as the identity on $F_{\lambda_i}(T)$, we have for any $x\in F_{\lambda_i}(T)$ and $y\in \bigoplus_{j\neq i}F_{\lambda_j}(T)$,
$$
\norm{x+p_n(\lambda_i\inv T)(y)}=\norm{p_n(\lambda_i\inv T)(x+y)} \leqslant \norm{x+y}.
$$
Note also that
$$
(\Id-\lambda_i\inv T)p_n(\lambda_i\inv T)=n^{-1}(\Id-\lambda_i^{-n}T^n).
$$
Since $T-\Id$ is injective on $\bigoplus_{j\neq i}F_{\lambda_j}(T)$ and therefore also invertible on $\bigoplus_{j\neq i}F_{\lambda_j}(T)$, we have for $y\in \bigoplus_{j\neq i}F_{\lambda_j}(T)$,
$$
p_n(\lambda_i\inv T)(y)=n^{-1}\big[(\Id-\lambda_i\inv T)|_{\bigoplus_{j\neq i}F_{\lambda_j}(T)}\big]\inv(\Id-\lambda_i^{-n}T^n)(y),
$$
whereby since $(\lambda_i^{-n} T^n)_n$ is bounded, $\lim_{n \rightarrow \infty}p_n(\lambda_i\inv T)(y)=0$. Applying this to the inequality $\norm{x+p_n(\lambda_i\inv T)(y)} \leqslant \norm{x+y}$, we get
$\|P_i(x+y)\|=\|x\| \leqslant \|x+y\|$.
\end{proof}

\begin{lemme}\label{technical}
Suppose $X$ is a complex space. Then for any $S,T\in {\rm Isom}_f(X)$, $\lambda\in \sigma(T)$ and $\mu_1,\ldots,\mu_k\in \sigma(S)$, if
$$
F_\lambda(T)\cap \bigoplus_{i=1}^kF_{\mu_i}(S)\neq\{0\},
$$
then
$$
\min_i|\mu_i-\lambda|\leqslant k\|S-T\|.
$$
It follows that for any $S,T\in {\rm Isom}_f(X)$ and $\lambda\in \sigma(T)$,
$$
\dim F_\lambda(T)\;\leqslant\; \dim\bigoplus_{|\mu-\lambda|\leqslant |\sigma(S)|\cdot\|S-T\|}F_\mu(S).
$$
\end{lemme}

\begin{proof}
Suppose that $x_i\in F_{\mu_i}(S)$ and $x_1+\ldots+x_k\in F_\lambda(T)$.
Then, by Lemma \ref{norm projection}, we have
\[\begin{split}
\min_i|\mu_i-\lambda|\cdot\|x_1+\ldots+x_k\|
&\leqslant \min_i|\mu_i-\lambda|\big(\|x_1\|+\ldots+\|x_k\|\big)\\
&\leqslant k\, \min_i|\mu_i-\lambda|\cdot\max_j\|x_j\|\\
&\leqslant k\, \max_j\|(\mu_j-\lambda)x_j\|\\
&\leqslant k\, \|(\mu_1-\lambda)x_1+\ldots+(\mu_k-\lambda)x_k\|\\
&=k\,\|S(x_1+\ldots+x_k)-T(x_1+\ldots+x_k)\|\\
&\leqslant k\,\|S-T\|\cdot \|x_1+\ldots+x_k\|.
\end{split}\]
Dividing by $\|x_1+\ldots+x_k\|$, we get $\min_i|\mu_i-\lambda|\leqslant k\|S-T\|$.

Now, suppose towards a contradiction that $A=\big\{\mu\in \sigma(S)\del |\mu-\lambda|\leqslant |\sigma(S)|\cdot\|S-T\|\big\}$ and
$$
\dim F_\lambda(T)\;>\; \dim\bigoplus_{\mu\in A}F_\mu(S).
$$
Then, as $X=\bigoplus_{\mu\in \sigma(S)}F_\mu(S)$, we must have
$$
F_\lambda(T)\,\cap\, \bigoplus_{\mu \in \sigma(S)\setminus A}F_\mu(S)\neq \{0\}
$$
and so $|\mu-\lambda|\leqslant  |\sigma(S)|\cdot\|S-T\|$ for some $\mu\in \sigma(S)\setminus A$, which is absurd.
\end{proof}

\begin{lemme}\label{continuity of spectra}
Let $\ku K(\T)$ denote the set of compact subsets of $\T$ equipped with the Vietoris topology, i.e., the topology induced by the Hausdorff metric, and let $X$ be a complex space. Then $T\mapsto \sigma(T)$ is norm continuous as a map from ${\rm Isom}_f(X)$ to $\ku K(\T)$.
\end{lemme}

\begin{proof}
Using that  $GL(X)$ is a norm open subset of $\ku L(X)$, we see that if $S_n\Lim{n\til \infty} T$ and $\lambda_n\Lim{n\til \infty}\lambda$ for $\lambda_n\in \sigma(S_n)$, then also $\lambda\in \sigma(T)$.
It follows that if $\ku U\subseteq \T$ is open and $\sigma(T)\subseteq \ku U$, then $\sigma(S)\subseteq \ku U$ for all $S$ in a neighbourhood of $T$. So $T\mapsto \sigma(T)$ is at least semi-continuous.

For the other half of continuity, note first that by Lemma \ref{uniform continuity} there is some $N$ such that  $|\sigma(S)|\leqslant N$ for all $S$ in some neighbourhood $\ku N$ of $T$. Therefore, by Lemma \ref{technical},  for all $S\in \ku N$ and $\lambda\in \sigma(T)$,
$$
\min_{\mu \in \sigma(S)}|\mu-\lambda|\leqslant N\|S-T\|.
$$
It follows that if $\ku U\subseteq \T$ is open and $\sigma(T)\cap \ku U\neq\tom$, then also $\sigma(S)\cap \ku U\neq\tom$ for all $S$ in a neighbourhood of $T$, finishing the proof of continuity.
\end{proof}

Finally, several examples of non-reflexive HI spaces have been considered in the literature.  Gowers  \cite{Go} constructed a separable HI space such that every subspace has non-separable dual and  Argyros, A. Arvanitakis and A. Tolias \cite{Ar} contructed a non-separable, (necessarily non-reflexive) HI space. Also, the Argyros--Haydon space $AH$ \cite{AH} is HI and has dual isomorphic to $\ell_1$. As noted, any single isometry of a complex HI space $X$ must act nearly trivially on $X$.  But nothing is known about the global behaviour of ${\rm Isom}(X)$ even in the above mentioned examples.


\subsection{Problems concerning Hilbert space}
As mentioned in the introduction, the following strong versions of the second part of Mazur's rotation problem remain open.
\begin{prob}
Suppose $\norm\cdot$ is an equivalent maximal or almost transitive norm on $\ku H$. Must $\norm\cdot$  be euclidean?
\end{prob}

Also, though not every bounded subgroup $G\leqslant GL(\ku H)$ permits an equivalent $G$-invariant euclidean norm, the following is open.
\begin{prob}
Suppose $G\leqslant GL(\ku H)$ is a bounded subgroup. Must there be an equivalent $G$-invariant maximal, almost transitive or even transitive norm on $\ku H$?
\end{prob}

Note that by Theorem \ref{conv trans} convex transitivity coincides with almost transitivity on $\ku H$.

In  \cite{CS} it is mentioned that if $X$ is a space with an almost transitive norm and there exists a finite-dimensional isometry $\Id+F$ for which $F$ has rank $1$, then $X$ is isometric to a Hilbert space. However, the following question is still open.

\begin{prob}[Cabello-S\'anchez] Let $X$ be a space with an (almost) transitive norm, and which admits a non-trivial finite-dimensional isometry. Must $X$ be Hilbertian?
\end{prob}

Note, that by Theorem \ref{convextransitive}, if $X$ is separable reflexive and satisfies the hypothesis, then $X$ must have a Schauder basis.


\subsection{General questions}
Of course, to prove Conjecture \ref{conj HI}, one is tempted to eliminate the second option in Theorem \ref{decompos}. For this, the following would be a intermediate step.

\begin{prob} 
Let $X$ be a Banach space,  $G\leqslant {\rm Isom}_f(X)$ and assume that 
$$
\sup \{\dim F_T\del  T \in G\}<\infty.
$$
Must $G$ act nearly trivially on $X$?
\end{prob}

On the other hand, for a potential counter-example, one might begin with the following.
\begin{prob}
Find a separable space $X$ and a bounded subgroup of
$GL_f(X)$ which is infinite and discrete for SOT.
\end{prob}

We know of no real uniformly convex space for which no renorming is maximal. Considered as a real space, the example of Theorem \ref{nomaximal} is also HI and uniformly convex, but does not satisfy the $\lambda \Id+S$-property.

\begin{prob} 
Does there exist a real super-reflexive Banach space without a maximal norm? Without an almost transitive norm?
\end{prob}

As mentioned, the fact that there exist complex super-reflexive spaces with no equivalent almost transitive norms, shows that  a certain approach to smooth renormings does not work. Our counter-examples are therefore candidates for a negative answer to the following question.

\begin{prob}
 Does there exist a Banach space $X$ and constants $1 \leqslant q \leqslant 2\leqslant p$ such that the set of equivalent norms on $X$ with modulus of  convexity of type $p$ and the set of equivalent norms on $X$ with modulus of smoothness of type $q$ are disjoint and both non-empty? Is the space defined in \cite{F} or one of its subspace an example of this?
\end{prob}

Also, one may ask whether the two hypotheses of having a convex transitive norm and having a non-trivial finite-dimensional isometry in Theorem \ref{convextransitive} can be separated, so as dealing with different norms. I.e., we would be interested in the answer to the following.
\begin{prob}
Let $X$ be a separable, reflexive, Banach space with a convex transitive norm. Does it follow that $X$ has a Schauder basis?
\end{prob}


\end{document}